\documentclass{article}
\usepackage{amsfonts}
\usepackage{mathrsfs}
\usepackage{dsfont}
\usepackage{bbding}
\usepackage{amssymb}
\usepackage{amsmath}
\usepackage{graphicx}
\usepackage{booktabs}
\usepackage{multirow}
\usepackage{hyperref}
\usepackage{indentfirst}
\usepackage{subfigure}

\hypersetup{
hidelinks,
colorlinks=true,
linkcolor=blue,
citecolor=blue
}

\newcommand\bbf[1]{\mathchoice{\hbox{\boldmath$\displaystyle#1$}}
{\hbox{\boldmath$\textstyle#1$}} {\hbox{\boldmath$\scriptstyle#1$}}
{\hbox{\boldmath$\scriptscriptstyle#1$}} }
\newcommand{\bea}{\begin{eqnarray}}
\newcommand{\eea}{\end{eqnarray}}
\newcommand{\Bea}{\begin{eqnarray*}}
\newcommand{\Eea}{\end{eqnarray*}}
\newcommand{\ba}{\begin{array}}
\newcommand{\ea}{\end{array}}
\newcommand{\bt}{\begin{tabular}}
\newcommand{\et}{\end{tabular}}
\newcommand{\btb}{\begin{table}}
\newcommand{\etb}{\end{table}}
\newcommand{\bc}{\begin{center}}
\newcommand{\ec}{\end{center}}
\newcommand{\beq}{\begin{equation}}
\newcommand{\eeq}{\end{equation}}

\begin{document}

\title{Updatable Estimation in Generalized Linear Models with Missing Response
 }
\author{Xianhua Zhang\textsuperscript{a}, Lu Lin\textsuperscript{a,}\footnote{The corresponding
author. Email: linlu@sdu.edu.cn; sdulinlu@163.com. The research was
supported by National Key R\&D Program of China (2018YFA0703900) and NNSF project (11971265) of China.} , Qihua Wang\textsuperscript{b}
\\
\small \textsuperscript{a}Zhongtai Securities Institute for Financial Studies, Shandong University, Jinan, China
\\
\small \textsuperscript{b}Academy of Mathematics and Systems Sciences, Chinese Academy of Sciences, Beijing, China
}

\date{}
\maketitle

\begin{abstract}
\large\baselineskip=18pt
This paper develops an updatable inverse probability weighting (UIPW) estimation for the generalized linear models with response missing at random in streaming data sets.
A two-step online updating algorithm is provided for the proposed method.
In the first step we construct an updatable estimator for the parameter in propensity function and hence obtain an updatable estimator of the propensity function;
in the second step we propose an UIPW estimator with the inverse of the updating  propensity function value at each observation as the weight for estimating the parameter of interest.
The UIPW estimation is universally applicable due to its relaxation on the constraint on the number of data batches.
It is shown that the proposed estimator is consistent and asymptotically normal with the same asymptotic variance as that of the oracle estimator, and hence the oracle property is obtained.
The finite sample performance of the proposed estimator is illustrated by the simulation and real data analysis.
All numerical studies confirm that the UIPW estimator performs as well as the batch learner.
\end{abstract}
~\\
{\it \textbf{Keywords}}:\large{ Generalized linear model; Streaming data sets; Missing response; Updatable inverse probability weighting.}

\baselineskip=20pt
\newpage
\section{Introduction}
With the development of modern science and technology, massive data sets arise in various fields.
It not only presents opportunities, but also brings challenges, such as the difficulty of storing on standard computers and the insufficient computation speed of traditional algorithms.
Until now, the major strategies to meet challenges have been researched by many statisticians, most of which can be divided into the three categories:
sub-sampling (Liang et al., 2013; Kleiner et al., 2014; Ping et al., 2015), divide and conquer (Lin \& Xi, 2011; Scott et al., 2016; Chen \& Xie, 2014; Song \& Liang, 2014; Pillonetto et al., 2019) and the online updating (Schifano et al, 2016; Wang et al., 2018; Xue et al., 2020; Luo \& Song, 2020).
However, the online updating is basically different from the others because the streaming data set is a special data type in the emerging field of "big data".
More specifically, in the streaming data set, the data arrive in streams and chunks sequentially. So the statistical method should be of online updating framework, without storage requirement for historical data.

The stochastic gradient descent (SGD) algorithm proposed by Robbins \& Monro (1951) is a well-known iterative algorithm.
Since the SGD algorithm updates the estimation with one data point at each time, it coincides with the mechanism of how data are observed in streaming data sets and becomes a natural solution to the online estimation.
The SGD algorithm has been developed by statisticians with many improved versions, such as averaged implicit SGD (AISGD) (Toulis et al., 2014) and SGD-QN (Bordes et al., 2009).
The online second-order method is another popular approach for handling streaming data sets, such as Amari et al. (2000) proposed the online Newton step.
Furthermore, Chen et al. (2021) proposed an SGD-based online algorithm that can make decisions and update the decision rule online.
In addition to the SGD-based online algorithm, Schifano et al. (2016) extended the cumulative estimating equation (CEE) estimator (Lin \& Xi, 2011) and proposed the cumulatively updated estimating equation (CUEE) estimator.
The CUEE estimator had been shown to be superior to other divide-and-conquer or online-updated estimators in terms of bias and mean squared error.
However, the CUEE estimator sacrifices application value for the oracle property and estimation consistency.
For example, the total number of streaming data sets is usually constrained by $K=O\left(N^{c}_K\right)$, where $0\textless c \textless 1$ is a constant, $N_K$ is the number of total data and $K$ is the total number of the data batches.
This constraint means that the total number of the data batches cannot be very large.
Lin \& Zhang (2002) proposed a two-sample empirical Euclidean penalty likelihood method based on historical empirical estimation and current empirical information.
Luo \& Song (2020) proposed the renewable estimation (RE) through an incremental updating algorithm which can relax the constraint.
But the above methods only focused on the case of complete data.

As Zhu et al. (2019) pointed out, in the era of big data, it is more likely to encounter incomplete observations.
Imputation (Rubin, 1987; 1996) is a widely adopted approach for dealing with missing data.
Wang (2008) proposed two nonparametric probability density estimations: inverse probability weighted estimation and the regression calibration density estimation.
However, the estimating methods faced the challenge of curse of dimensionality.
Hu et al. (2010, 2012) considered dimension reduction in estimation of response mean to avoid the curse of dimensionality.
Based on the kernel smoothing and sieve method, Wang et al. (2021) developed two distributed nonparametric imputation methods, which extended the divide and conquer strategy to missing response problems.
However, applying these imputation methods directly to the streaming data set is inappropriate (Syavasya \& Lakshmi, 2022).
To handle the missing data in the streaming data set, Wellenzohn et al. (2017) proposed Top-k Case Matching (TKCM) model.
It consistently imputes the missing data and retains accuracy when there is a block of missing data.
Peng et al. (2019) developed Incremental Space-Time based Model (ISTM) to impute the missing data in the streaming data set.
Inverse probability weighting (IPW) (Horvitz \& Thompson, 1952; Robins et al., 1994) is another popular approach for handling missing data.
The inverse probability weighted method is applied to missing response problem to define M-estimators in Wang (2007).
However, the existing IPW approaches for handling missing data are of the offline framework.
How to use IPW approach to deal with the missing data in streaming data sets is a significant problem.
To the best of our knowledge, there is little relevant literature to solve it. 
It is desired to develop online updating IPW estimation approaches and algorithms.


In this paper, we propose an online updating IPW estimation approach for the generalized linear model (GLM) (McCullagh \& Nelder, 1983) in the presence of missing response.
The proposed approach is developed with the commonly used IPW approach. Our contributions include the followings:
\begin{itemize}
    \item[(a)] A two-step online updating approach is proposed such that both the inverse probability weights and the estimator of the parameter of interest can update simultaneously;
    \item[(b)] We establish the asymptotic property of the updatable estimator without the condition $K=O(N_K^{c}),0\textless c\textless 1$;
    \item[(c)] Extend the proposed approach to the heterogeneous streaming data set. We get the efficient score function by the projection method, and then use a two-step online updating approach to get the online updating estimation.
\end{itemize}

The remainder of this paper is organized as follows.
In Section 2, the UIPW estimation is proposed to construct the updatable estimator in the streaming data set with missing response for the GLM, via a two-step online updating algorithm.
The theoretical properties of the updatable estimator are investigated.
We extend the UIPW estimation to the heterogeneous streaming data set in Section 3 with an efficient UIPW (EUIPW) estimation.
The main simulation studies and real data analysis are provided in Section 4 to evaluate the proposed method.
Technical details are provided in Appendix.

\section{Updatable inverse probability weighting estimation}
\subsection{Methodology}
Let $W=(Y,X^{\top})^{\top}$ denote the full variable vector and $\boldsymbol{\beta} \in \boldsymbol{\Theta}_{\boldsymbol{\beta}}$ denote the $p$-dimensional parameter vector of interest with true value $\boldsymbol{\beta}^0$, where $Y$ is the response, $X \in \mathbb{R}^{p}$ is the associated covariate.
In the complete case, for the GLM with $g(E[Y|X])=X^{\top}\boldsymbol{\beta}$, we assume that all the data $W_i=(Y_i,X_i^{\top})^{\top}$ for $i=1,2,\ldots,n$ are independent and identically distributed and follow exponential family distribution with the probability density function (PDF)
\begin{align}
    f\left(y ; \eta, \phi\right)=\exp \left\{\frac{y \eta-b\left(\eta\right)}{\phi}+c\left(y, \phi\right)\right\},\label{eq:pdf}
\end{align}
where $\phi$ is the dispersion parameter.
Suppose that $\eta=\eta(x^{\top}\boldsymbol{\beta})$ is a monotone differentiable function.
In the GLM, we have the following properties:
\begin{align}
    \begin{split}&
        \frac{\partial b(\eta(x^{\top}\boldsymbol{\beta}))}{\partial \eta(x^{\top}\boldsymbol{\beta})} = \mu,\text{ } \frac{\partial \mu}{\partial \eta(x^{\top}\boldsymbol{\beta})} = v(\mu),
        \\&
        \frac{\partial \eta(x^{\top}\boldsymbol{\beta})}{\partial \boldsymbol{\beta}}=\frac{\partial \eta(x^{\top}\boldsymbol{\beta})}{\partial \mu}\frac{\partial \mu}{\partial g(\mu)}\frac{\partial g(\mu)}{\partial \boldsymbol{\beta}}=[v(\mu)g_{\mu}(\mu)]^{-1}x^{\top},
    \end{split}\label{eq:GLM properties}
\end{align}
where $\mu=E[Y|X=x]$, $v(\mu)$ is the known unit variance function defined as $Var(Y|X=x)=\phi  v\left(\mu\right)$ and $g_{\mu}(\cdot)$ is the derivative of $g(\cdot)$ with respect to $\mu$.
We write the associated log-likelihood function as:
$$
\ell\left(\boldsymbol{\beta}, \phi ; \bbf W_n\right)=\frac{1}{n}\sum_{i=1}^n \log f\left(Y_{i} ; X_{i}, \boldsymbol{\beta}, \phi\right)=\frac{1}{n}\sum_{i=1}^n \left\{\frac{Y_{i} \eta(X_{i}^{\top} \boldsymbol{\beta})-b\left(\eta(X_{i}^{\top} \boldsymbol{\beta})\right)}{\phi }+c\left(Y_{i}, \phi\right)\right\},
$$
where $\bbf W_n=(W_1,\ldots,W_n)$.
Maximizing $\ell\left(\boldsymbol{\beta}, \phi ; \bbf W_n\right)$ is equivalent to maximizing $(1/n)\sum_{i=1}^n \{Y_{i} \eta(X_{i}^{\top} \boldsymbol{\beta})-b\left(\eta(X_{i}^{\top} \boldsymbol{\beta})\right)\}$ in $\ell\left(\boldsymbol{\beta}, \phi ; \bbf W_n\right)$ with respect to $\boldsymbol{\beta}$.
According to (\ref{eq:GLM properties}), the corresponding score function is $(1/n)\sum_{i=1}^n(Y_i-\mu_i)[v(\mu_i)g_{\mu}(\mu_i)]^{-1}X_i^{\top}$.

In the presence of missing data, we assume the response $Y$ is subject to missingness, but the covariate $X$ is fully observed.
Let $\delta$ denote the indicator of observing $Y$, i.e., $\delta=1$ if $Y$ is observed and $\delta=0$ otherwise.
The observed data are $n$ independent and identically distributed copies of $(\delta, Y, X^{\top})$.
We assume the response $Y$ is missing at random (MAR) (Rubin, 1976), that is
\begin{equation}
    P(\delta=1 | Y, X)=P(\delta=1 | X), \label{eq:propensity function}
\end{equation}
which means that the response is missing depending on covariates only but not on itself. 

Let $\pi(x)=P(\delta=1|X=x)$.
However, in most situations, we need to deduce $\pi(X)$ based on the data.
In particular, we generally can posit a model for $P(\delta=1|X)$; for example, a full parametric model $\pi(X;\boldsymbol{\alpha})$, where $\pi(\cdot;\cdot)$ is a known function, but $\boldsymbol{\alpha} \in \boldsymbol{\Theta}_{\boldsymbol{\alpha}}$ is an unknown parameter vector.
With $\boldsymbol{\alpha}$ estimated by $\hat{\boldsymbol{\alpha}}$, the IPW estimator for $\boldsymbol{\beta}$ can be obtained according to Horvitz \& Thompson (1952).

We consider the case where the samples arrive sequentially in chunks and the response is missing at random.
With slightly abusing the notation, let $W=(\delta,Y,X^{\top})^{\top}$.
The series of $k$ batches of data are denoted by $W_1=(\delta_1,Y_1,X_1^{\top})^{\top},\ldots,$\\$W_{N_k}=(\delta_{N_k},Y_{N_k},X_{N_k}^{\top})^{\top}$, where $W_i$ are i.i.d. copies of $W$, totally $N_k$ observations available.
For $j=1,\ldots,k,\ldots$, let $\boldsymbol{d}_j=(W_i:i\in \bbf i_j)$ be sequential data sets with the index sets $\bbf i_j$ defined by $\bbf i_j=\{N_{j-1}+1,\ldots,N_{j-1}+n_j\}$ with $N_0=0$ and $N_{k-1}+n_k=N_k$.
When the $k$-th batch of data arrives, consider the IPW approach, for any $j\in\{1,\ldots,k\}$, the score function on the $j$-th batch of data is
$$
\mathbf{S}(\boldsymbol{\beta}|\bbf d_j,\boldsymbol{\alpha})=\frac{1}{n_j}\sum_{i \in \bbf i_j}S(W_i;\boldsymbol{\beta},\boldsymbol{\alpha})=\frac{1}{n_j}\sum_{i \in \bbf i_j}\frac{\delta_i}{\pi(X_i;\boldsymbol{\alpha})}(Y_i-\mu_i)[v(\mu_i)g_{\mu}(\mu_i)]^{-1}X_i^{\top},
$$

Our goal is to construct the updatable estimator of $\boldsymbol{\beta}$ under the above settings.
To this end, the first task is to obtain the online updating estimation of the unknown parameter vector $\bbf \alpha$ in the propensity function $\pi(X;\boldsymbol{\alpha})$. 
Note that according to the binomial distribution, the log-likelihood function for $\bbf \alpha$ on the $j$-th batch of data for any $j \in \{1,\ldots,k\}$ is 
$$
\frac{1}{n_j}\sum_{i \in \bbf i_j} \log\bigg\{\pi(X_i;\boldsymbol{\alpha})^{\delta_i}\big[1-\pi(X_i;\boldsymbol{\alpha})\big]^{(1-\delta_i)}\bigg\}.
$$ 
Thus, the corresponding estimating function on the $j$-th batch of data $\boldsymbol{d}_j$ is
\begin{equation}
    \mathbf{V}(\boldsymbol{\delta}_j;\mathbf{X}_j,\boldsymbol{\alpha})=\frac{1}{n_j}\sum_{i \in \bbf i_j}V(\delta_i;X_i,\boldsymbol{\alpha})=\frac{1}{n_j}\sum_{i \in \bbf i_j} \frac{\nabla_{\boldsymbol{\alpha}}\pi(X_i;\boldsymbol{\alpha})\big[\delta_i-\pi(X_i;\boldsymbol{\alpha})\big]}{\pi(X_i;\boldsymbol{\alpha})\big[1-\pi(X_i;\boldsymbol{\alpha})\big]} = 0, \label{eq:alpha estimation equation}
\end{equation}
where $\boldsymbol{\delta}_j=\{\delta_i,i \in \bbf i_j\}$, $\mathbf{X}_j=\{X_i,i \in \bbf i_j\}$ and $\nabla_{\boldsymbol{\alpha}}\pi(X_i;\boldsymbol{\alpha})$ is the derivative of $\pi(X_i;\boldsymbol{\alpha})$ with respect to $\bbf \alpha$.

The remainder task is to construct the updatable estimator of $\bbf\beta$.
Let $\mathbf{R}_{\boldsymbol{\beta}}(\boldsymbol{\beta}|d,\boldsymbol{\alpha})=\nabla_{\boldsymbol{\beta}}\mathbf{S}(\boldsymbol{\beta}|d,\boldsymbol{\alpha})$ denote the derivative of $\mathbf{S}(\boldsymbol{\beta}|d,\boldsymbol{\alpha})$ with respect to $\boldsymbol{\beta}$, $\mathbf{R}_{\boldsymbol{\alpha}}(\boldsymbol{\beta}|d,\boldsymbol{\alpha})=\nabla_{\boldsymbol{\alpha}}\mathbf{S}(\boldsymbol{\beta}|d,\boldsymbol{\alpha})$ denote the derivative of $\mathbf{S}(\boldsymbol{\beta}|d,\boldsymbol{\alpha})$ with respect to $\boldsymbol{\alpha}$ and $\mathbf{R}_{\boldsymbol{\alpha}\boldsymbol{\beta}}(\boldsymbol{\beta}|d,\boldsymbol{\alpha})=\nabla_{\boldsymbol{\alpha}}\nabla_{\boldsymbol{\beta}}\mathbf{S}(\boldsymbol{\beta}|d,\boldsymbol{\alpha})$ stand for the second-order mixed derivatives of $\mathbf{S}(\boldsymbol{\beta}|d,\boldsymbol{\alpha})$ with respect to $\boldsymbol{\beta}$ and $\boldsymbol{\alpha}$.
We use the IPW estimator of $\boldsymbol{\beta}$ on the whole data as the oracle estimator.
When the second batch of data $\boldsymbol{d}_2$ arrives after the first batch of data $\boldsymbol{d}_1$, we need to update the initial IPW estimator $\hat{\boldsymbol{\beta}}_{N_1}^*=\hat{\boldsymbol{\beta}}_1$ to a renewed IPW estimator $\hat{\boldsymbol{\beta}}_{N_2}^*$, without using of the raw data in $\boldsymbol{d}_1$.
Here, the IPW estimation $\hat{\boldsymbol{\beta}}_1$ satisfies the score equation $\mathbf{S}(\boldsymbol{\beta}|\bbf d_1,\hat{\boldsymbol{\alpha}}_1)=0$.
Suppose we have obtained REs $\hat{\boldsymbol{\alpha}}_1$ and $\hat{\boldsymbol{\alpha}}_2$.
In the GLM, IPW estimator $\hat{\boldsymbol{\beta}}_{N_2}^*$ is the solution to the equation 
\begin{align}
    \mathbf{S}(\boldsymbol{\beta}|\bbf d_1,\hat{\boldsymbol{\alpha}}_2)+\mathbf{S}(\boldsymbol{\beta}|\bbf d_2,\hat{\boldsymbol{\alpha}}_2)=0. \label{eq:k=2}
\end{align}
However, after the raw data set $\bbf d_1$ is discarded, the above equation cannot be employed directly. 
Similar to Lin \& Zhang (2002) and Luo \& Song (2020), we introduce the following online updating strategy.
Taking the first-order Taylor series expansion of $\mathbf{S}(\hat{\boldsymbol{\beta}}|\bbf d_1,\hat{\boldsymbol{\alpha}}_2)$ around $(\hat{\boldsymbol{\beta}}_1,\hat{\boldsymbol{\alpha}}_1)$ in (\ref{eq:k=2}), we propose a new estimator $\hat{\boldsymbol{\beta}}_2$ as the solution to the equation
\begin{align}
    \mathbf{R}_{\boldsymbol{\alpha}}(\hat{\boldsymbol{\beta}}_1|\bbf d_1,\hat{\boldsymbol{\alpha}_1})(\hat{\boldsymbol{\alpha}}_2-\hat{\boldsymbol{\alpha}}_1)+\big[\mathbf{R}_{\boldsymbol{\beta}}(\hat{\boldsymbol{\beta}}_1|\bbf d_1,\hat{\boldsymbol{\alpha}}_1)+\mathbf{R}_{\boldsymbol{\alpha}\boldsymbol{\beta}}(\hat{\boldsymbol{\beta}}_1|\bbf d_1,\hat{\boldsymbol{\alpha}}_1)(\hat{\boldsymbol{\alpha}}_2-\hat{\boldsymbol{\alpha}}_1)\big]\nonumber\\
    (\boldsymbol{\beta}-\hat{\boldsymbol{\beta}}_1)+\mathbf{S}(\boldsymbol{\beta}|\bbf d_2,\hat{\boldsymbol{\alpha}}_2)=0. \label{eq:k=2 estimating equation}
\end{align}
In the above, the statistics $\hat{\boldsymbol{\alpha}}_1$, $\hat{\boldsymbol{\beta}}_1$, $\mathbf{R}_{\boldsymbol{\alpha}}(\hat{\boldsymbol{\beta}}_1|\bbf d_1,\hat{\boldsymbol{\alpha}_1})$, $\mathbf{R}_{\boldsymbol{\beta}}(\hat{\boldsymbol{\beta}}_1|\bbf d_1,\hat{\boldsymbol{\alpha}_1})$ and $\mathbf{R}_{\boldsymbol{\alpha}\boldsymbol{\beta}}(\hat{\boldsymbol{\beta}}_1|\bbf d_1,\hat{\boldsymbol{\alpha}_1})$ were stored before $\bbf d_2$ arrives.
We then use the stored statistics and $\bbf d_2$ to solve the equation (\ref{eq:k=2 estimating equation}), and the resulting solution is defined to be the estimator of $\boldsymbol{\beta}$, $\hat{\boldsymbol{\beta}}_2$ say, after $\bbf d_2$ arrives.
The solution to equation (\ref{eq:k=2 estimating equation}) can be obtained by the Newton-Raphson algorithm.

Generalizing the above procedure to streaming data sets, we can construct the updatable estimator of $\boldsymbol{\beta}$ by the following estimating equation:
\begin{align}
   \sum_{j=1}^{k-1}\mathbf{L}_{1}(\hat{\boldsymbol{\beta}}_j|\bbf d_j,\hat{\boldsymbol{\alpha}}_j)(\hat{\boldsymbol{\alpha}}_k-\hat{\boldsymbol{\alpha}}_{k-1})+\sum_{j=1}^{k-1}\mathbf{L}_{2}(\hat{\boldsymbol{\beta}}_j|\bbf d_j,\hat{\boldsymbol{\alpha}}_j)(\boldsymbol{\beta}-\hat{\boldsymbol{\beta}}_{k-1})+\mathbf{S}(\boldsymbol{\beta}|\bbf d_k,\hat{\boldsymbol{\alpha}}_k)=0, \label{eq:update beta}
\end{align}
where
\begin{align*}
    \begin{split}&
        \mathbf{L}_{1}(\hat{\boldsymbol{\beta}}_j|\bbf d_j,\hat{\boldsymbol{\alpha}}_j)=\mathbf{R}_{\boldsymbol{\alpha}}(\hat{\boldsymbol{\beta}}_j|\bbf d_j,\hat{\boldsymbol{\alpha}}_j)+\mathbf{R}_{\boldsymbol{\alpha}\boldsymbol{\beta}}(\hat{\boldsymbol{\beta}}_j|\bbf d_j,\hat{\boldsymbol{\alpha}}_j)(\hat{\boldsymbol{\beta}}_{k-1}-\hat{\boldsymbol{\beta}}_j),\\
        &\mathbf{L}_{2}(\hat{\boldsymbol{\beta}}_j|\bbf d_j,\hat{\boldsymbol{\alpha}}_j)=\mathbf{R}_{\boldsymbol{\beta}}(\hat{\boldsymbol{\beta}}_j|\bbf d_j,\hat{\boldsymbol{\alpha}}_j)+\mathbf{R}_{\boldsymbol{\alpha}\boldsymbol{\beta}}(\hat{\boldsymbol{\beta}}_j|\bbf d_j,\hat{\boldsymbol{\alpha}}_j)(\hat{\boldsymbol{\alpha}}_{k}-\hat{\boldsymbol{\alpha}}_j).    
    \end{split}
\end{align*}
Solving equation (\ref{eq:update beta}) can be easily done by the following incremental
updating algorithm:
\begin{align}
    \hat{\boldsymbol{\beta}}_{k}^{(r+1)}=\hat{\boldsymbol{\beta}}_{k}^{(r)}-\bigg\{\tilde{\mathbf{L}}_{k-1}+\mathbf{R}_{\boldsymbol{\beta}}(\hat{\boldsymbol{\beta}}_{k-1}|\bbf d_k,\hat{\boldsymbol{\alpha}}_k)\bigg\}^{-1}\tilde{\mathbf{F}}_k^{(r)}, \label{eq:beta incremental}
\end{align}
where
\begin{equation*}\begin{split}&
    \tilde{\mathbf{L}}_k= \sum_{j=1}^{k}\mathbf{L}_{2}(\hat{\boldsymbol{\beta}}_j|\bbf d_j,\hat{\boldsymbol{\alpha}}_j),
\\&
    \tilde{\mathbf{F}}_{k}^{(r)}=\sum_{j=1}^{k-1}\mathbf{L}_{1}(\hat{\boldsymbol{\beta}}_j|\bbf d_j,\hat{\boldsymbol{\alpha}}_j)(\hat{\boldsymbol{\alpha}}_k-\hat{\boldsymbol{\alpha}}_{k-1})+\sum_{j=1}^{k-1}\mathbf{L}_{2}(\hat{\boldsymbol{\beta}}_j|\bbf d_j,\hat{\boldsymbol{\alpha}}_j)(\hat{\boldsymbol{\beta}}_{k}^{(r)}-\hat{\boldsymbol{\beta}}_{k-1})+\mathbf{S}(\hat{\boldsymbol{\beta}}_{k}^{(r)}|\bbf d_k,\hat{\boldsymbol{\alpha}}_k).
\end{split}\end{equation*}

Actually, the above is a two-step online updating algorithm. In the first step we obtain the updatable estimator $\hat{\boldsymbol{\alpha}}_k$ for the
parameter vector $\boldsymbol{\alpha}_k$; in the second step,
the updatable estimator $\hat{\boldsymbol{\beta}}_k$ can then be constructed with $\hat{\boldsymbol{\alpha}}_k$. The procedure of constructing the updatable estimator is summarized as follows:

{\it Step 1}: Choose an initial value of $\hat{\boldsymbol{\alpha}}_0 = \bbf 0$ and suppose that the first $(k-1)$-th online updating estimators $\hat{\boldsymbol{\alpha}}_1,\ldots, \hat{\boldsymbol{\alpha}}_{k-1}$ are obtained.
Then, the $k$-th online updating estimator $\hat{\boldsymbol{\alpha}}_k$  can be attained as  the solution to the following estimating equation:
\begin{equation}
    \sum_{j=1}^{k-1} \mathbf{H}\left(\boldsymbol{\delta}_j;\mathbf{X}_{j},\hat{\boldsymbol{\alpha}}_{j} \right) \left(\boldsymbol{\alpha}-\hat{\boldsymbol{\alpha}}_{k-1}\right)-\mathbf{V}\left(\boldsymbol{\delta}_k;\mathbf{X}_{k},\boldsymbol{\alpha} \right)=0, \label{eq:update alpha}
\end{equation} where $\mathbf{H}(\boldsymbol{\delta}_j;\mathbf{X}_j,\boldsymbol{\alpha})=-\nabla_{\boldsymbol{\alpha}}\mathbf{V}(\boldsymbol{\delta}_j;\mathbf{X}_j,\boldsymbol{\alpha})$.
The incremental updating algorithm for the above equation will be given in (\ref{eq:alpha incremental}).

{\it Step 2}:
Choose an initial value of $\hat{\boldsymbol{\beta}}_0 = \bbf 0$ and suppose that the first $(k-1)$-th updatable estimators $\hat{\boldsymbol{\beta}}_1,\ldots, \hat{\boldsymbol{\beta}}_{k-1}$ are obtained.
Then, by $\hat{\boldsymbol{\alpha}}_k$, $\hat{\boldsymbol{\alpha}}_{k-1}$, $\hat{\boldsymbol{\beta}}_{k-1}$, $\sum_{j=1}^{k-2}\mathbf{L}_{1}(\hat{\boldsymbol{\beta}}_j|\bbf d_j,\hat{\boldsymbol{\alpha}}_j)$ and $\sum_{j=1}^{k-2}\mathbf{L}_{2}(\hat{\boldsymbol{\beta}}_j|\bbf d_j,\hat{\boldsymbol{\alpha}}_j)$, the updatable estimator $\hat{\boldsymbol{\beta}}_k$  can be defined as  the solution to (\ref{eq:update beta}), and we can calculate the updatable estimator $\hat{\boldsymbol{\beta}}_k$ by (\ref{eq:beta incremental}).

We repeat the above steps until a stopping rule is met.
Note that solving equation (\ref{eq:update alpha}) may be easily done by the following incremental updating algorithm:
\begin{equation}
    \hat{\boldsymbol{\alpha}}_{k}^{(r+1)}=\hat{\boldsymbol{\alpha}}_{k}^{(r)}+\left\{\tilde{\mathbf{H}}_{k-1}+\mathbf{H}\left(\boldsymbol{\delta}_k;\mathbf{X}_{k},\hat{\boldsymbol{\alpha}}_{k-1}\right)\right\}^{-1} \tilde{\mathbf{V}}_{k}^{(r)}, \label{eq:alpha incremental}
\end{equation}
where $\tilde{\mathbf{H}}_{k}=\sum_{j=1}^{k} \mathbf{H}\left(\boldsymbol{\delta}_j;\mathbf{X}_{j},\hat{\boldsymbol{\alpha}}_{j}\right)$ and $\tilde{\mathbf{V}}_{k}^{(r)}=\tilde{\mathbf{H}}_{k-1}\left(\hat{\boldsymbol{\alpha}}_{k-1}-\hat{\boldsymbol{\alpha}}_{k}^{(r)}\right)+\mathbf{V}\left(\boldsymbol{\delta}_k;\mathbf{X}_{k},\hat{\boldsymbol{\alpha}}_{k}^{(r)}\right)$.
The above is an online updating form because it only involves the current data $\boldsymbol{d}_k$, the previous estimators $\hat{\boldsymbol{\alpha}}_{k-1}$ and $\hat{\boldsymbol{\beta}}_{k-1}$ together with the accumulative quantity $\sum_{j=1}^{k-1} \mathbf{H}\left(\boldsymbol{\delta}_j;\mathbf{X}_{j},\hat{\boldsymbol{\alpha}}_{j}\right)$, $\sum_{j=1}^{k-1}\mathbf{L}_{1}(\hat{\boldsymbol{\beta}}_j|\bbf d_j,\hat{\boldsymbol{\alpha}}_j)$ and $\sum_{j=1}^{k-1}\mathbf{L}_{2}(\hat{\boldsymbol{\beta}}_j|\bbf d_j,\hat{\boldsymbol{\alpha}}_j)$.

\subsection{Theoretical properties}
In this section, we establish the consistency and asymptotic normality for the proposed updatable estimators under the homogeneous model in Section 2.1, and then show its asymptotic equivalence to the oracle estimator.

For an arbitrary batch $k$, suppose that $W_1,\ldots,W_{N_k}$ are $i.i.d.$ samples from exponential family distribution with density $f(y ; x, \boldsymbol{\beta}, \phi)$, $i = 1,\ldots,N_k$.
Due to the MAR mechanism, given $X$, random variables $\delta$ and $Y$ are independent.
Let $\boldsymbol{\alpha}^0$ be the true value of the parameter $\boldsymbol{\alpha}$. The Fisher information matrix for $\boldsymbol{\alpha}$ is
\begin{align*}
    \mathcal{I}_{N_{k}}\left(\boldsymbol{\alpha}\right)=\sum_{j=1}^{k} E\left\{\mathbf{V}^{\top} (\boldsymbol{\delta}_j;\mathbf{X}_j,\boldsymbol{\alpha})\mathbf{V}(\boldsymbol{\delta}_j;\mathbf{X}_j,\boldsymbol{\alpha})\right\}=\sum_{i=1}^{N_k} E\big\{X_i\pi(X_i;\boldsymbol{\alpha})[1-\pi(X_i;\boldsymbol{\alpha})]X_i^{\top}\big\},
\end{align*}
and the Fisher information matrix for $\boldsymbol{\beta}$ is
\begin{equation*}\begin{split}&
    \mathcal{I}_{N_{k}}\left(\boldsymbol{\beta}\right)=\sum_{j=1}^{k} E\left\{\mathbf{S}^{\top} (\boldsymbol{\beta}|\bbf d_j,\boldsymbol{\alpha}^0) \mathbf{S}(\boldsymbol{\beta}|\bbf d_j,\boldsymbol{\alpha}^0)\right\}=\sum_{i=1}^{N_k}E\big\{X_i [\pi(X_i;\boldsymbol{\alpha}^0)v\left(\mu_i\right)g^2_{\mu}(\mu_i)]^{-1}X_i^{\top}\big\}.
\end{split}
\end{equation*}

We assume the following regularity conditions for establishing the asymptotic properties of $\hat{\boldsymbol{\alpha}}_k$:\\
C1. For a fixed positive constant $c$, there exists an open subset $\omega_{\boldsymbol{\alpha}}$ of $\Theta_{\boldsymbol{\alpha}}$ containing the true parameter point $\boldsymbol{\alpha}^0$ such that $\inf_x\pi(x;\boldsymbol{\alpha})\textgreater c \textgreater 0$ for all $\boldsymbol{\alpha} \in \omega_{\boldsymbol{\alpha}}$.\\
C2. $\mathcal{I}_{N_k}(\boldsymbol{\alpha})$ is finite and positive definite for all $\boldsymbol{\alpha} \in \omega_{\boldsymbol{\alpha}}$.\\
C3. $\pi(x;\boldsymbol{\alpha})$ is twice continuously differentiable and $\pi(x;\boldsymbol{\alpha})$ is Lipschitz continuous in $\omega_{\boldsymbol{\alpha}}$.

For establishing the asymptotic properties of $\hat{\boldsymbol{\beta}}_k$, we introduce the following regularity conditions:\\
C4. $\sup_x Var(Y|X=x)\textless \infty$.\\
C5. There exists an open subset $\omega_{\boldsymbol{\beta}}$ of $\Theta_{\boldsymbol{\beta}}$ containing the true parameter point $\boldsymbol{\beta}^0$ such that $\mathcal{I}_{N_k}(\boldsymbol{\beta})$ is positive definite for all $\boldsymbol{\beta} \in \omega_{\boldsymbol{\beta}}$.\\
C6. There exist functions $M_{stl}$ such that
$$
\left|\frac{\partial^3}{\partial \beta_s \partial \beta_t \partial \beta_l} b(\eta(x^{\top}\boldsymbol{\beta}))\right| \leq M_{s t l}(x)\text { for all } \boldsymbol{\beta} \in \omega_{\boldsymbol{\beta}},
$$
where
$$
E_{\boldsymbol{\beta}^0}\left[M_{s t l}(X)\right]<\infty \text { for all } s, t, l.
$$
C7. $\ell(\boldsymbol{\beta};d)$ is twice continuously differentiable with respect to $\boldsymbol{\beta}$.
$\mathbf{R}_{\boldsymbol{\beta}}(\boldsymbol{\beta}|d,\boldsymbol{\alpha})$ is Lipschitz continuous in $\omega_{\boldsymbol{\beta}}$ and $\mathbf{R}_{\boldsymbol{\alpha}}(\boldsymbol{\beta}|d,\boldsymbol{\alpha})$ is Lipschitz continuous in $\omega_{\boldsymbol{\alpha}}$.\\
C8. $\sup_{\boldsymbol{\beta}\in \omega_{\boldsymbol{\beta}}}E\left\{X\left[\pi(X;\boldsymbol{\alpha})v(\mu)g_{\mu}^2(\mu)\right]^{-1}X^{\top}\right\}<\infty$.

Condition C1 is necessary for the consistency of the estimator $\hat{\boldsymbol{\alpha}}_k$ and ensures the finite asymptotic variance of the estimator $\hat{\boldsymbol{\alpha}}_k$.
Conditions C2, C4, C5 and C6 are the standard regularity conditions (see, e.g., Lehmann EL, 1983).
Conditions C3 and C7 are needed for renewable estimators, see Luo \& Song (2020).
Condition C8 is needed to establish the asymptotic normality of the UIPW estimator.

\textbf{Lemma 1.} Under conditions C1-C3, the updatable estimator $\hat{\boldsymbol{\alpha}}_k$ given in equation (\ref{eq:update alpha}) is asymptotically normally distributed, i.e.
\begin{align}
  \sqrt{N_k}(\hat{\boldsymbol{\alpha}}_k-\boldsymbol{\alpha}^0) \overset{d}{\rightarrow} N\bigg(0,\bigg\{E\big[V(\delta;X,\boldsymbol{\alpha}^0)V(\delta;X,\boldsymbol{\alpha}^0)^{\mathrm{T}}\big]\bigg\}^{-1}\bigg)\text{ as }N_k=\sum_{j=1}^{k}n_j \rightarrow \infty.\nonumber
\end{align}

For proof of Lemma 1 see Luo \& Song (2020).

\textbf{Theorem 2.} Under conditions C1-C7, the updatable estimator $\hat{\boldsymbol{\beta}}_k$ given in equation (\ref{eq:update beta}) is consistent, namely, $\hat{\boldsymbol{\beta}}_k \overset{P}{\rightarrow} \boldsymbol{\beta}^0$ as $N_k=\sum _{j=1}^kn_j\rightarrow \infty$. 

The proof of Theorem 2 is given in Appendix A.1.

\textbf{Theorem 3.} Under conditions C1-C8, the updatable estimator $\hat{\boldsymbol{\beta}}_k$ is asymptotically normal, i.e.
\begin{align}
    \sqrt{N_k}(\hat{\boldsymbol{\beta}}_k-\boldsymbol{\beta}^0) \overset{d}{\rightarrow} N(0,\boldsymbol{\Sigma}^0)\text{ as }N_k=\sum_{j=1}^{k}n_j \rightarrow \infty,\label{eq:beta asymptotically normal}
\end{align}
where
\begin{align*}
    &\boldsymbol{\Sigma}^0=\bigg\{E\left[X\left(v(\mu)g_{\mu}^2(\mu)\right)^{-1}X^{\top}\right]\bigg\}^{-1}
    Var\bigg\{S(\boldsymbol{\beta}^0|W,\boldsymbol{\alpha}^0)
    -J(\boldsymbol{\beta}^0,\boldsymbol{\alpha}^0)\bigg\}
    \\
    &\qquad\enspace\bigg\{E\left[X\left(v(\mu)g_{\mu}^2(\mu)\right)^{-1}X^{\top}\right]^{\top}\bigg\}^{-1},\\
    &J(\boldsymbol{\beta}^0,\boldsymbol{\alpha}^0)
    =E\big[S(\boldsymbol{\beta}^0|W,\boldsymbol{\alpha}^0)
    V(\delta;X,\boldsymbol{\alpha}^0)^{\top}\big]\bigg\{E\big[V(\delta;X,\boldsymbol{\alpha}^0)V(\delta;X,\boldsymbol{\alpha}^0)^{\top}\big]\bigg\}^{-1}V(\delta;X,\boldsymbol{\alpha}^0).
\end{align*}

The proof of Theorem 3 is provided in Appendix A.2. 

\textbf{Remark 1.} 
From the Theorem 2, it can be seen that the updatable estimator has the standard convergence rate of order $O_p(N_k^{-1/2})$, and the asymptotic covariance is the same as that of the oracle estimator.
Moreover, the theoretical result always holds without any constraint on $k$. It is an important property for the GLM with streaming data sets because it means that the method is adaptive to the situation where streaming data sets arrive fast and perpetually.

It is attractive that the asymptotic covariance matrix of UIPW estimator $\hat{\boldsymbol{\beta}}_k$ is the same as the oracle estimator $\hat{\boldsymbol{\beta}}_{N_k}^*$.
This implies that we construct a fully efficient estimator using some summary statistics that are convenient to store, without the need to use historical data in the process, which greatly relaxes the constraint on the memory of computers.

\textbf{Remark 2.}
In many cases, the asymptotic covariance matrix $\boldsymbol{\Sigma}_0$ needs to be estimated, such as constructing confidence intervals.
Therefore, we estimate $\boldsymbol{\Sigma}_0$ by replacing the terms under expectation in (\ref{eq:beta asymptotically normal}) with the empirical averages
evaluated at $(\hat{\boldsymbol{\beta}}_k,\hat{\boldsymbol{\alpha}}_k)$, expressed as
\begin{align*}
&\hat{\boldsymbol{\Sigma}}_k=
    \bigg\{\widehat{E}\left[X\left(v(\mu)g_{\mu}^2(\mu)\right)^{-1}X^{\top}\right]\bigg\}^{-1}\widehat{Var}\bigg\{S
    (\hat{\boldsymbol{\beta}}_k|W,\hat{\boldsymbol{\alpha}}_k)-
    \hat{J}(\hat{\boldsymbol{\beta}}_k,\hat{\boldsymbol{\alpha}}_k)\bigg\}
    \\
    &\qquad\enspace\bigg\{\widehat{E}\left[X\left(v(\mu)g_{\mu}^2(\mu)\right)^{-1}X^{\top}\right]^{\top}\bigg\}^{-1},\\
&\hat{J}(\hat{\boldsymbol{\beta}}_k,\hat{\boldsymbol{\alpha}}_k)
=\hat{E}\big[S(\hat{\boldsymbol{\beta}}_k|W,\hat{\boldsymbol{\alpha}}_k)
V(\delta;X,\hat{\boldsymbol{\alpha}}_k)^{\top}\big]\bigg\{\hat{E}\big[V(\delta;X,\hat{\boldsymbol{\alpha}}_k)V(\delta;X,\hat{\boldsymbol{\alpha}}_k)^{\top}\big]\bigg\}^{-1}\mathbf{V}(\delta;X,\hat{\boldsymbol{\alpha}}_k),
\end{align*} 
where $\widehat{E}[\cdot]$ is the empirical averages, for instance,
$$
\widehat{E}\left[X\left(v(\mu)g_{\mu}^2(\mu)\right)^{-1}X^{\top}\right]=\frac{1}{n_k}\sum_{i\in \bbf i_k}\left\{X_i\left[v(\mu_i)g_{\mu}^2(\mu_i)\right]^{-1}X_i^{\top}\right\}.
$$
Then, $\hat{\boldsymbol{\Sigma}}_k$ is the consistent estimator of $\boldsymbol{\Sigma}^0$.

Based on the asymptotic distribution of the UIPW estimator in Theorem 2, the Wald test can be used to test hypotheses of coefficients.
Define the null hypothesis as $\mathbf{h}_0:\boldsymbol{\beta}=\boldsymbol{\beta}^{\mbox{null}}$, where $\boldsymbol{\beta}^{\mbox{null}} \in \Theta_{\boldsymbol{\beta}}$ is a pre-fixed null $p$-dimension vector.
Under the null hypothesis, the Wald test statistic is
\begin{align*}
    (\hat{\boldsymbol{\beta}}_k-\boldsymbol{\beta}^{\mbox{null}})^{\top} \hat{\boldsymbol{\Sigma}}_k^{-1}(\hat{\boldsymbol{\beta}}_k-\boldsymbol{\beta}^{\mbox{null}}) \overset{d}{\rightarrow} \chi_p^2,
\end{align*}
where $\chi_p^2$ is a Chi-square distribution with $p$ degrees of freedom.
Thus, if the confidence level is $1-c$, a $100(1-c)\%$ confidence ellipsoid for $\boldsymbol{\beta}$ is given by
$$
\mathcal{C}=\{\boldsymbol{\beta}:(\hat{\boldsymbol{\beta}}_k-\boldsymbol{\beta})^{\top} \hat{\boldsymbol{\Sigma}}_k^{-1}(\hat{\boldsymbol{\beta}}_k-\boldsymbol{\beta})<\chi_p^2(\alpha)\}.
$$

\section{Extension to the case of heterogeneous distributions}
The aforementioned algorithm needs the condition that the data in different data sets are identically distributed.
However, the heterogeneity of distribution across different batches of data in a streaming data set is a common situation.
Usually, we can use a heterogeneous parameter vector $\bbf\gamma_j$ to characterize the heterogeneity (Duan et al., 2021).
Thus, in the GLM, for $j \in \{1,\ldots,k,\ldots\}$, we write the associated log-likelihood function as $\ell(\boldsymbol{\beta},\boldsymbol{\gamma}_j;\bbf d_j)$.
The score functions with respect to $\boldsymbol{\beta}$ and $\boldsymbol{\gamma}_j$ can be expressed as $\mathbf{S}(\boldsymbol{\beta}|\boldsymbol{d}_j,\boldsymbol{\alpha},\boldsymbol{\gamma}_j)$ and $\mathbf{T}(\boldsymbol{\gamma}_j|\boldsymbol{d}_j,\boldsymbol{\beta},\boldsymbol{\alpha})=(1/n_j)\sum_{i \in \bbf i_j}T(W_i;\boldsymbol{\beta},\boldsymbol{\alpha},\boldsymbol{\gamma}_j)$, respectively.
Note that the true value of $\boldsymbol{\gamma}_j$ may be different for each batch of data in the streaming data set. 
In below we present two examples, both of which have an appreciable degree of generality.

\textit{Example 1.}
Suppose that the true value of dispersion parameter $\phi$ in (\ref{eq:pdf}) is different among different data blocks.
For $j \in \{1,\ldots,k,\ldots\}$, let $\phi_j=\gamma_j$ and $\mathbf{I}_j$ to be the $n_j \times n_j$ identity matrix, the samples are independent and identically distributed and follow the exponential family distribution with the PDF
$$
f\left(y ; \eta, \gamma_j\right)=\exp \left\{\frac{y \eta-b\left(\eta\right)}{\gamma_j}+c\left(y, \gamma_j\right)\right\},
$$
where $\eta=\eta(x^{\top}\boldsymbol{\beta})$. The GLM is
$$
Y=g^{-1}(X^{\top}\boldsymbol{\beta})+\boldsymbol{\epsilon}_j,
$$
where $Var[\boldsymbol{\epsilon}_j]=\gamma_jv(\mu)\mathbf{I}_j$.
Then the score functions for the common parameter vector $\boldsymbol{\beta}$ and the nuisance parameter $\boldsymbol{\gamma}_j$ at the $j$-th batch of data are
\begin{align*}\begin{split}&
    \mathbf{S}(\boldsymbol{\beta}|\boldsymbol{d}_j,\boldsymbol{\alpha},\boldsymbol{\gamma}_j)=\frac{1}{n_j}\sum_{i \in \bbf i_j}S(W_i;\boldsymbol{\beta},\boldsymbol{\alpha},\boldsymbol{\gamma}_j)=\frac{1}{n_j\gamma_j}\sum_{i \in \bbf i_j}\frac{\delta_i}{\pi(X_i;\boldsymbol{\alpha})}(Y_i-\mu_i)[v_j(\mu_i)g_{\mu}(\mu_i)]^{-1}X_i^{\top},\\
    &\mathbf{T}(\boldsymbol{\gamma}_j|\boldsymbol{d}_j,\boldsymbol{\beta},\boldsymbol{\alpha})=\frac{1}{n_j}\sum_{i \in \bbf i_j}T(W_i;\boldsymbol{\beta},\boldsymbol{\alpha},\boldsymbol{\gamma}_j)=\frac{1}{n_j}\sum_{i \in \bbf i_j}\bigg[-\frac{Y_i\eta(X_i^{\top}\boldsymbol{\beta})-b(\eta(X_i^{\top}\boldsymbol{\beta}))}{\boldsymbol{\gamma}_j^{2}}+\frac{\partial c(Y_i,\boldsymbol{\gamma_j})}{\partial \boldsymbol{\gamma_j}}\bigg],
\end{split}
\end{align*}
respectively.

\textit{Example 2.} 
In this example, suppose that the full data vector is $W=(\delta,Y,X^{\top},Z^{\top})^{\top}$, where $\delta$ is the indicator variable, Y is the response, X and Z are associated covariates. 
We still assume that the response is missing depending on $X$.
For $j \in \{1,\ldots,k,\ldots\}$, we consider the samples are independent and identically distributed and follow the exponential family distribution with the PDF
\begin{align}
    f\left(y ; \eta_j, \phi\right)=\exp \left\{\frac{y \eta_j-b\left(\eta_j\right)}{\phi}+c\left(y, \phi\right)\right\}.\label{eq:example2 pdf}
\end{align}
The difference between (\ref{eq:pdf}) and(\ref{eq:example2 pdf}) is that we suppose $\eta_j=\eta(x^{\top}\boldsymbol{\beta}+z^{\top}\boldsymbol{\gamma}_j)$ in (\ref{eq:example2 pdf}).
The GLM is
$$
g(E[Y|X,Z])=X^{\top}\boldsymbol{\beta}+Z^{\top}\boldsymbol{\gamma}_j.
$$
Similar to Section 2.1, for $j \in \{1,\ldots,k,\ldots\}$, we have $\mu=E[Y|X=x,Z=z]$ and $Var[Y|X=x,Z=z]=\phi v(\mu)$.
The score functions for the common parameter vector $\boldsymbol{\beta}$ and the nuisance parameter $\boldsymbol{\gamma}_j$ at the $j$-th batch of data are
\begin{align*}\begin{split}&
    \mathbf{S}(\boldsymbol{\beta}|\boldsymbol{d}_j,\boldsymbol{\alpha},\boldsymbol{\gamma}_j)=\frac{1}{n_j}\sum_{i \in \bbf i_j}S(W_i;\boldsymbol{\beta},\boldsymbol{\alpha},\boldsymbol{\gamma}_j)=\frac{1}{n_j}\sum_{i \in \bbf i_j}\frac{\delta_i}{\pi(X_i;\boldsymbol{\alpha})}(Y_i-\mu_i)[v(\mu_i)g_{\mu}(\mu_i)]^{-1}X_i^{\top},\\
    &\mathbf{T}(\boldsymbol{\gamma}_j|\boldsymbol{d}_j,\boldsymbol{\beta},\boldsymbol{\alpha})=\frac{1}{n_j}\sum_{i \in \bbf i_j}T(W_i;\boldsymbol{\beta},\boldsymbol{\alpha},\boldsymbol{\gamma}_j)=\frac{1}{n_j}\sum_{i \in \bbf i_j}\frac{\delta_i}{\pi(X_i;\boldsymbol{\alpha})}(Y_i-\mu_i)[v(\mu_i)g_{\mu}(\mu_i)]^{-1}Z_i^{\top},
\end{split}
\end{align*}
respectively.

Compared with Section 2.1, we add a nuisance parameter in score function to describe the heterogeneity such that the score function $S(W;\boldsymbol{\beta},\boldsymbol{\alpha},\boldsymbol{\gamma}_j)$ has two nuisance parameters.
Let $\boldsymbol{\theta}_j=(\boldsymbol{\beta},\boldsymbol{\alpha},\boldsymbol{\gamma}_j) \in \mathbb{R}^d$.
In this section, $\boldsymbol{\beta} \in \mathbb{R}^p$ is a common parameter vector of interest on all the batches of data, $\boldsymbol{\alpha}\in\mathbb{R}^p$ is a common nuisance parameter vector on all the batches of data, and the $(d-2p)$-dimensional nuisance parameter $\boldsymbol{\gamma}_j$ is allowed to be different across batches of data. The true value of $\boldsymbol{\theta}_j$ is denoted by $\boldsymbol{\theta}^0_{j}$.
To deal with the heterogeneity, we propose the following efficient updatable inverse probability weighting (EUIPW) to estimate the common parameters.

Since $\boldsymbol{\gamma}_j$ is the nuisance parameter, we propose to approximate the efficient score function. 
Motivated by theories of efficient score, the efficient score function is defined as the projection of the score function of $\boldsymbol{\beta}$ on the space that is orthogonal to the space spanned by the score function of nuisance parameter $\boldsymbol{\gamma}_j$ (van der Vaart, 1998).
In our setting, it can be expressed as
\begin{align}
    \mathbf{U}(\boldsymbol{d}_j|\boldsymbol{\beta},\boldsymbol{\alpha},\boldsymbol{\gamma}_j) = \mathbf{S}(\boldsymbol{\beta}|\boldsymbol{d}_j,\boldsymbol{\alpha},\boldsymbol{\gamma}_j)-\mathcal{I}_{\boldsymbol{\beta}\boldsymbol{\gamma}}^{(j)} \mathcal{I}_{\boldsymbol{\gamma}\boldsymbol{\gamma}}^{(j)} \mathbf{T}(\boldsymbol{\gamma}_j|\boldsymbol{d}_{j},\boldsymbol{\beta},\boldsymbol{\alpha}), \label{eq:ES}
\end{align}
for $j \in \{1,\ldots,k,\ldots\}$, where $\mathcal{I}_{\boldsymbol{\beta}\boldsymbol{\gamma}}^{(j)}$ and $\mathcal{I}_{\boldsymbol{\gamma}\boldsymbol{\gamma}}^{(j)}$ are the corresponding submatrices of the information matrix for the $j$-th batch of data which are expressed as
\begin{align}\begin{split}&
    \mathcal{I}_{\boldsymbol{\beta}\boldsymbol{\gamma}}^{(j)} = E\bigg\{\mathbf{S}(\boldsymbol{\beta}|\boldsymbol{d}_j,\boldsymbol{\alpha},\boldsymbol{\gamma}_j)\mathbf{T}^{\top}(\boldsymbol{\gamma}_j|\boldsymbol{d}_j,\boldsymbol{\beta},\boldsymbol{\alpha})\bigg\},\\
    &\mathcal{I}_{\boldsymbol{\gamma}\boldsymbol{\gamma}}^{(j)} = E\bigg\{\mathbf{T}(\boldsymbol{\gamma}_j|\boldsymbol{d}_j,\boldsymbol{\beta},\boldsymbol{\alpha})\mathbf{T}^{\top}(\boldsymbol{\gamma}_j|\boldsymbol{d}_j,\boldsymbol{\beta},\boldsymbol{\alpha})\bigg\}. \nonumber
\end{split}
\end{align}

When the $k$-th batch of data arrives, we only update the estimators of $\boldsymbol{\beta}$ and $\boldsymbol{\alpha}$.
The estimator of $\boldsymbol{\gamma}_k$ can be obtained based on the $k$-th batch of data. The online updating procedure is designed as follows:

{\it Step 1}: As in Step 1 in Section 2, we can obtain $\hat{\boldsymbol{\alpha}}_{k-1}$ and $\hat{\boldsymbol{\alpha}}_k$, then the estimation of $\boldsymbol{\gamma}_k$ can be derived by solving the following equation:
\begin{align*}
    \mathbf{S}(\boldsymbol{\beta}|\boldsymbol{d}_k,\hat{\boldsymbol{\alpha}}_k,\boldsymbol{\gamma}_k) = 0,
    \mathbf{T}(\boldsymbol{\gamma}_k|\boldsymbol{d}_{k},\boldsymbol{\beta},\hat{\boldsymbol{\alpha}}_k) = 0 \mbox{ for } \bbf\beta \mbox{ and } \bbf\gamma_k.
\end{align*}
Note that the estimation of $\boldsymbol{\beta}$ obtained in this step is not used in the following.

{\it Step 2}: Suppose that the first $(k-1)$-th updatable estimators $\hat{\boldsymbol{\beta}}_1,\ldots, \hat{\boldsymbol{\beta}}_{k-1}$  are obtained. Then, by $\hat{\boldsymbol{\alpha}}_k$, $\hat{\boldsymbol{\alpha}}_{k-1}$, $\hat{\boldsymbol{\beta}}_{k-1}$, $\hat{\boldsymbol{\gamma}}_k$ together with $\sum_{j=1}^{k-2}\mathbf{G}_{1}(\hat{\boldsymbol{\beta}}_j|\bbf d_j,\hat{\boldsymbol{\alpha}}_j,\hat{\boldsymbol{\gamma}}_j)$ and $\sum_{j=1}^{k-2}\mathbf{G}_{2}(\hat{\boldsymbol{\beta}}_j|\bbf d_j,\hat{\boldsymbol{\alpha}}_j,\hat{\boldsymbol{\gamma}}_j)$, the updatable estimator $\hat{\boldsymbol{\beta}}_k$  can be attained as  the solution to the following estimating equation:
\begin{align}
    \sum_{j=1}^{k-1}\mathbf{G}_{1}(\hat{\boldsymbol{\beta}}_j|\bbf d_j,\hat{\boldsymbol{\alpha}}_j,\hat{\boldsymbol{\gamma}}_j)(\hat{\boldsymbol{\alpha}}_k-\hat{\boldsymbol{\alpha}}_{k-1})+\sum_{j=1}^{k-1}\mathbf{G}_{2}(\hat{\boldsymbol{\beta}}_j|\bbf d_j,\hat{\boldsymbol{\alpha}}_j,\hat{\boldsymbol{\gamma}}_j)(\boldsymbol{\beta}-\hat{\boldsymbol{\beta}}_{k-1})\nonumber\\
    +\mathbf{U}(\boldsymbol{\beta}|\bbf d_k,\hat{\boldsymbol{\alpha}}_k,\hat{\boldsymbol{\gamma}}_k)=0, \label{eq:update EUIPW}
\end{align}
where
\begin{align*}
    \mathbf{G}_{1}(\hat{\boldsymbol{\beta}}_j|\bbf d_j,\hat{\boldsymbol{\alpha}}_j,\hat{\boldsymbol{\gamma}}_j)=\nabla_{\boldsymbol{\alpha}}\mathbf{U}(\hat{\boldsymbol{\beta}}_j|\bbf d_j,\hat{\boldsymbol{\alpha}}_j,\hat{\boldsymbol{\gamma}}_j)+\nabla_{\boldsymbol{\alpha}}\nabla_{\boldsymbol{\beta}}\mathbf{U}(\hat{\boldsymbol{\beta}}_j|\bbf d_j,\hat{\boldsymbol{\alpha}}_j,\hat{\boldsymbol{\gamma}}_j)(\hat{\boldsymbol{\beta}}_{k-1}-\hat{\boldsymbol{\beta}}_j),\\
    \mathbf{G}_{2}(\hat{\boldsymbol{\beta}}_j|\bbf d_j,\hat{\boldsymbol{\alpha}}_j,\hat{\boldsymbol{\gamma}}_j)=\nabla_{\boldsymbol{\beta}}\mathbf{U}(\hat{\boldsymbol{\beta}}_j|\bbf d_j,\hat{\boldsymbol{\alpha}}_j,\hat{\boldsymbol{\gamma}}_j)+\nabla_{\boldsymbol{\alpha}}\nabla_{\boldsymbol{\beta}}\mathbf{U}(\hat{\boldsymbol{\beta}}_j|\bbf d_j,\hat{\boldsymbol{\alpha}}_j,\hat{\boldsymbol{\gamma}}_j)(\hat{\boldsymbol{\alpha}}_{k}-\hat{\boldsymbol{\alpha}}_j).
\end{align*}

We repeat the above steps until a stopping rule is met. Similarly to equation (\ref{eq:beta incremental}),  solving equation (\ref{eq:update EUIPW}) can be easily done by the Newton-Raphson algorithm.

The above is online updating form because it only involves the current data $\boldsymbol{d}_k$, the previous estimators $\hat{\boldsymbol{\alpha}}_{k-1}$ and $\hat{\boldsymbol{\beta}}_{k-1}$ together with the accumulative quantity
$\sum_{j=1}^{k-1} \mathbf{H}\left(\mathbf{X}_{j};\hat{\boldsymbol{\alpha}}_{j}\right)$, $\sum_{j=1}^{k-1}\mathbf{G}_{1}(\hat{\boldsymbol{\beta}}_j|\bbf d_j,\hat{\boldsymbol{\alpha}}_j)$ and $\sum_{j=1}^{k-1}\mathbf{G}_{2}(\hat{\boldsymbol{\beta}}_j|\bbf d_j,\hat{\boldsymbol{\alpha}}_j)$.

\section{Numerical analysis}
\subsection{Empirical evidences}
In this section, simulation experiments are conducted to evaluate the proposed online updating approach under the following two typical models: the linear regression model and logistic regression model.
The model of propensity function $\pi(X)=P(\delta=1 \mid X)$ is chosen as the following logistic regression model:
\begin{align}
    \pi(X;\boldsymbol{\alpha})=\frac{1}{1+e^{-X^{\top}\boldsymbol{\alpha}}}, \label{eq:simulation propensity function}
\end{align}
where $\boldsymbol{\alpha}$ is an unknown parameter vector.
We compare the proposed UIPW estimator with the following one benchmark and two competing estimators:
\begin{itemize}
    \item The oracle estimator $\hat{\boldsymbol{\beta}}_{N_k}^*$: an estimator of $\boldsymbol{\beta}$ obtained by the whole data and the offline IPW method, which is employed as a benchmark.
    \item The simple average IPW estimator $\hat{\boldsymbol{\beta}}^{ave}$: let $\hat{\boldsymbol{\beta}}_j^{ave}$ be the IPW estimator of $\boldsymbol{\beta}$ computed on the $j$-th batch of data for $j=1,2,\ldots,K$, we then get the simple average of them as $\hat{\boldsymbol{\beta}}^{ave}=\frac{1}{K}\Sigma_{j=1}^{K}\hat{\boldsymbol{\beta}}_j^{ave}$;
    \item The naive IPW estimator $\hat{\boldsymbol{\beta}}^c$ constructed only on the current data.
\end{itemize}

We use the oracle estimator $\hat{\boldsymbol{\beta}}_{N_k}^*$ as the gold standard in all comparisons.
The estimation performance is measured with the mean squared error (MSE) and computation time (C.Time(s)).
As we all know, the choice of incremental algorithm affects the computation time.
In our simulation, we choose stochastic gradient descent (SGD) algorithm and Newton–Raphson algorithm to deal with different data settings.
The main advantages and disadvantages of SGD and Newton-Raphson algorithms are shown in Table \ref{incremental}.
In general, the convergence rate of Newton-Raphson algorithm is relatively fast.
However, when the sample size is small, the Newton–Raphson algorithm may perform badly due to inappropriate initial value setting.
Therefore, we must sacrifice convergence rate to guarantee the convergence.
\begin{table}[htbp]
    \centering
    \caption{The advantages and disadvantages of SGD algorithm and Newton-Raphson algorithm}
    \resizebox{\textwidth}{!}{
    \begin{tabular}{ccc}
        \toprule
        Algorithm&Advantages&Disadvantages\\
        \midrule
        SGD&$\cdot$Independent of the initial value &$\cdot$First order convergence\\
        \hline
        \multirow{2}*{Newton–Raphson}&\multirow{2}*{$\cdot$Second order convergence}&$\cdot$Dependence on initial values may result in non-convergence\\
        & &$\cdot$Complex computation \\
        \bottomrule
    \end{tabular}}
    \label{incremental}
\end{table}

In our simulation studies, the choice of incremental algorithm based on the following criteria:
\begin{itemize}
    \item[(1)] The Newton–Raphson algorithm is chosen in the process of online updating estimation.
    \item[(2)] The SGD algorithm is employed for calculating the simple average IPW estimate and the naive IPW estimate if the size of batch of data is small.
\end{itemize}

In each model, we compare the four estimators under two different scenarios.
In Scenario 1, we fix $N_K$ but change $K$ while in Scenario 2 we fix $K$ (or $n_k$) but change $n_k$ (or $K$).
The simulation results under Scenario 1 are presented in the following Tables while Scenario 2 are Figures.

\subsubsection{Homoscedastic linear regression model}
We write the $j$-th batch of data as $\bbf d_j=\{\boldsymbol{\delta}_j,\mathbf{Y}_j,\mathbf{X}_j\}$, $j = 1, \ldots, K$, where $\boldsymbol{\delta}_j=(\delta_1^{(j)},\ldots,\delta_{n_j}^{(j)})^{\top}$ is the indicator vector, $\mathbf{Y}_j=(\mathrm{y}_{1}^{(j)},\ldots,\mathrm{y}_{n_j}^{(j)})^{\top}$ is the response subject to missing at random and $\mathbf{X}_j=(\mathrm{x}_{1}^{(j)},\ldots,\mathrm{x}_{n_j}^{(j)})^{\top}$ is the $p$-dimension covariate fully observed, where $n_j$ is the sample size of $j$-th batch of data and $\mathrm{y}_{i}^{(j)} \mid \mathrm{x}_{i}^{(j)}$ are independently sampled from a Gaussian distribution.
Then, the streaming data set is in the form $\boldsymbol{d}_1, \boldsymbol{d}_2, \ldots, \boldsymbol{d}_k, \ldots$, with $\boldsymbol{d}_j=\{W_i,i \in \bbf i_j\}$.
We let $\bbf O$ be the index set of all the observed $Y$.
The score function and corresponding negative Hessian matrix for $j$-th batch of data $\boldsymbol{d}_j$ are respectively denoted by $\mathbf{S}(\boldsymbol{\beta}|\bbf d_j,\boldsymbol{\alpha})=(1/n_j)\mathbf{X}_j^{\top}\mathbf{M}_j(\boldsymbol{\alpha})(\mathbf{Y}_j-\mathbf{X}_j\boldsymbol{\beta})$ and $\mathbf{R}_{\boldsymbol{\beta}}(\boldsymbol{\beta}|\bbf d_j,\boldsymbol{\alpha})=(1/n_j)\mathbf{X}_j^{\top}\mathbf{M}_j(\boldsymbol{\alpha})\mathbf{X}_j$, where $\mathbf{M}_j(\boldsymbol{\alpha})=\mbox{diag}(\delta_i^{(j)}/(\pi(\mathrm{x}_{i}^{(j)};\boldsymbol{\alpha})),i=1,\ldots,n_j$.

In this example, when the $j$-th batch of data arrives, we consider the following linear regression model:
\begin{align}
    \mathrm{y}_i^{(j)}=2\mathrm{x}_{i1}^{(j)}+1.5\mathrm{x}_{i2}^{(j)}+\mathrm{x}_{i3}^{(j)}+0.5\mathrm{x}_{i4}^{(j)}+\epsilon_i^{(j)}, \quad  i \in \bbf O, \label{eq:linear model}
\end{align}
where $\mathrm{x}_{\cdot t}^{(j)} \sim U(0,1), t=1,2$, $\mathrm{x}_{\cdot t}^{(j)} \sim N(0,1), t=3,4$, $\epsilon^{(j)} \sim N(0,1)$, $\boldsymbol{\beta}^0=(2,1.5,1,0.5)$. The propensity function is set to be:
\begin{align*}
    \pi(\mathbf{X}_j;\boldsymbol{\alpha}^0)=\frac{1}{1+\exp(-(0.5\mathrm{x}_{i1}^{(j)}+\mathrm{x}_{i2}^{(j)}+1.5\mathrm{x}_{i3}^{(j)}+0.5\mathrm{x}_{i4}^{(j)}))},\quad i=1,\ldots,n_j,
\end{align*}
where $\boldsymbol{\alpha}^0=(0.5,1,1.5,0.5)$.

\begin{figure}[htbp]
	\caption{The MSE and CTime in linear regression model (\ref{eq:linear model}). Note that the MSE of the naive estimator is not depicted in (a), (b), (c) and (d) because it is out of the designed figures.}
	\centering
    \vspace{-0.2cm}
    \subfigbottomskip=-2cm
	\subfigure{
        \centering
        \includegraphics[width=12cm,height=6.5cm]{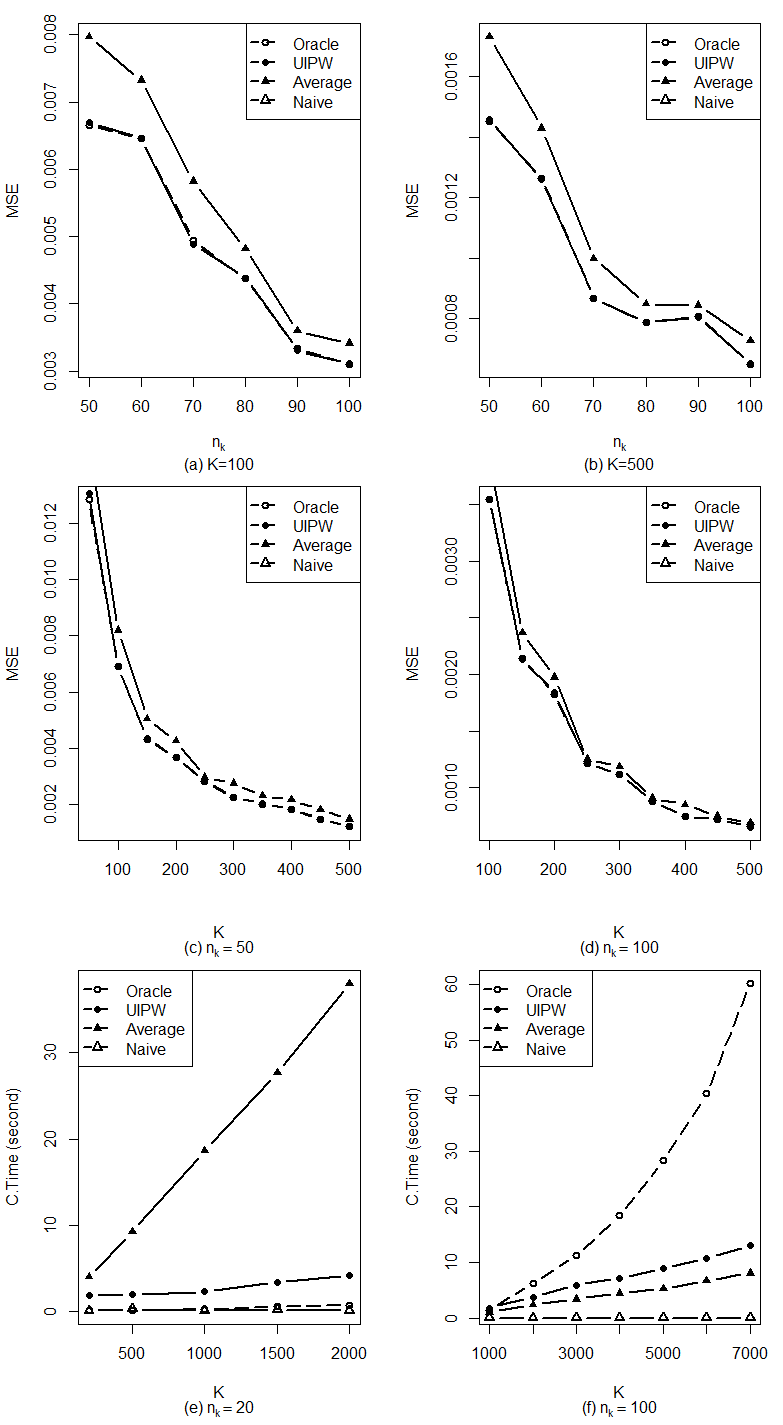}
	}

	\subfigure{
        \centering
		\includegraphics[width=12cm,height=6.5cm]{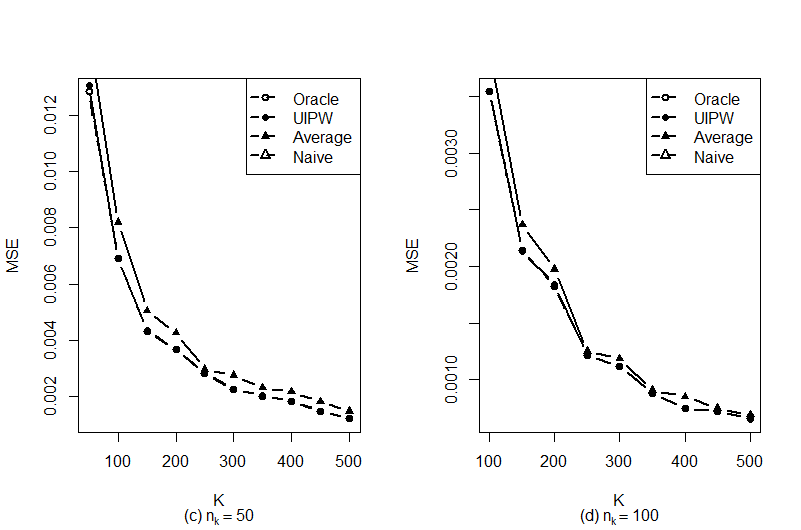}
	}
    \subfigure{
            \centering
			\includegraphics[width=12cm,height=6.5cm]{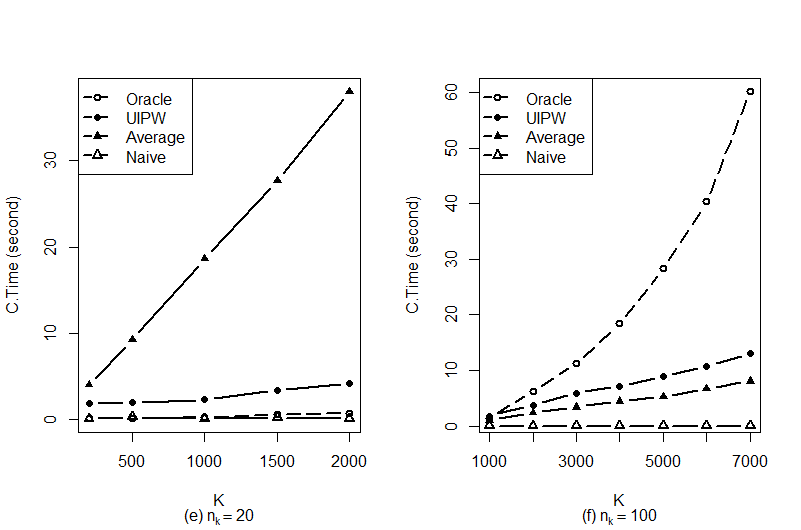}
	}
    \label{linearMSEp}
\end{figure}

The simulation results are reported in Table \ref{linearMSEt} and Figure \ref{linearMSEp}.
We have the following findings:
\begin{itemize}
    \item[(1)] \textit{Mean squared error.} Under the criterion of MSE,
    in Scenario 1, when $N_K$ is fixed, Table \ref{linearMSEt} indicates that the MSE of our UIPW estimator and the simple average IPW estimator increases as $K$ increases but our UIPW estimator appears fairly robust to different $K$.
    In Scenario 2, When $n_k$ or $K$ increases, as shown in Figure \ref{linearMSEp} (a)-(d), the MSE of our UIPW estimator decreases.
    In general, our UIPW estimator always exhibits similar performances to the oracle estimator and slightly better than the simple average IPW estimator. Unsurprisingly, the naive IPW estimator has the worst performance.
    \item[(2)] \textit{Computation time.} 
    Under the criterion of C.Time, the naive estimator takes the least C.Time because it only uses the current data.
    When $n_k$ is small, as shown in Figure \ref{linearMSEp} (e), the simple average IPW estimator takes the most C.Time, which increases rapidly as $K$ increases.
    Our UIPW estimator takes slightly more C.Time than the oracle estimator.
    However, setting a larger $n_k$, as shown in Figure \ref{linearMSEp} (f), the oracle estimator takes the most C.Time with $K$ increases, but our UIPW estimator still performs well.
\end{itemize}

\begin{table}[hbp]
    \centering
    \caption{The MSE and computation time of four estimators are summarized over 200 replication, under the setting of $N_K=100,000$ and $p = 4$ for the linear regression model (\ref{eq:linear model}) with batch size $n_k$ from 50 to 2000.}
    \begin{tabular}{|ccccc|}
        \hline
        \multirow{2}*{}         &\multicolumn{4}{c|}{$K=50$, $n_k=2000$}\\
                 &Oracle&UIPW&Average&Naive\\
        \hline
        MSE      &2.976$\times$$10^{-4}$&3.038$\times$$10^{-4}$&3.260$\times$$10^{-4}$&5.933$\times$$10^{-3}$\\
        C.Time(s)&0.2&11.7&8.08&0.17\\
        \hline
        \multirow{2}*{} &\multicolumn{4}{c|}{$K=100$, $n_k=1000$}\\
        &Oracle&UIPW&Average&Naive\\
        \hline
        MSE      &3.008$\times$$10^{-4}$&3.004$\times$$10^{-4}$&3.056$\times$$10^{-4}$&3.708$\times$$10^{-2}$\\
        C.Time(s)&0.25&6.1&4.59&0.04\\
        \hline
        \multirow{2}*{}         &\multicolumn{4}{c|}{$K=200$, $n_k=500$}\\
                 &Oracle&UIPW&Average&Naive\\
        \hline
        MSE      &3.272$\times$$10^{-4}$&3.271$\times$$10^{-4}$&3.402$\times$$10^{-4}$&6.588$\times$$10^{-2}$\\
        C.Time(s)&0.51&3.14&2.54&0.02\\
        \hline
        \multirow{2}*{} &\multicolumn{4}{c|}{$K=500$, $n_k=200$}\\
        &Oracle&UIPW&Average&Naive\\
        \hline
        MSE      &3.283$\times$$10^{-4}$&3.289$\times$$10^{-4}$&3.465$\times$$10^{-4}$&1.786$\times$$10^{-2}$\\
        C.Time(s)&0.73&1.64&1.12&0.01\\
        \hline
        \multirow{2}*{}         &\multicolumn{4}{c|}{$K=1000$, $n_k=100$}\\
                 &Oracle&UIPW&Average&Naive\\
        \hline
        MSE      &3.319$\times$$10^{-4}$&3.398$\times$$10^{-4}$&3.711$\times$$10^{-4}$&2.244$\times$$10^{-1}$\\
        C.Time(s)&1.38&2.06&1.3&0.01\\
        \hline
        \multirow{2}*{}&\multicolumn{4}{c|}{$K=2000$, $n_k=50$}\\
        &Oracle&UIPW&Average&Naive\\
        \hline
        MSE      &3.489$\times$$10^{-4}$&3.503$\times$$10^{-4}$&3.728$\times$$10^{-4}$&7.079$\times$$10^{-1}$\\
        C.Time(s)&2.7&2.79&1.54&0.01\\
        \hline
    \end{tabular}
    \label{linearMSEt}
\end{table}

\newpage
\subsubsection{Homoscedastic logistic regression model}
\begin{figure}[htbp]
	\caption{The MSE and computation time in the logistic regression model (\ref{eq:logistic model}). Note that the MSE of the naive estimator is not depicted in (a), (b), (c) and (d) because it is out of the designed figures.}
	\centering
    \vspace{-0.2cm}
    \subfigbottomskip=-2cm
	\subfigure{
        \centering
        \includegraphics[width=12cm,height=6.5cm]{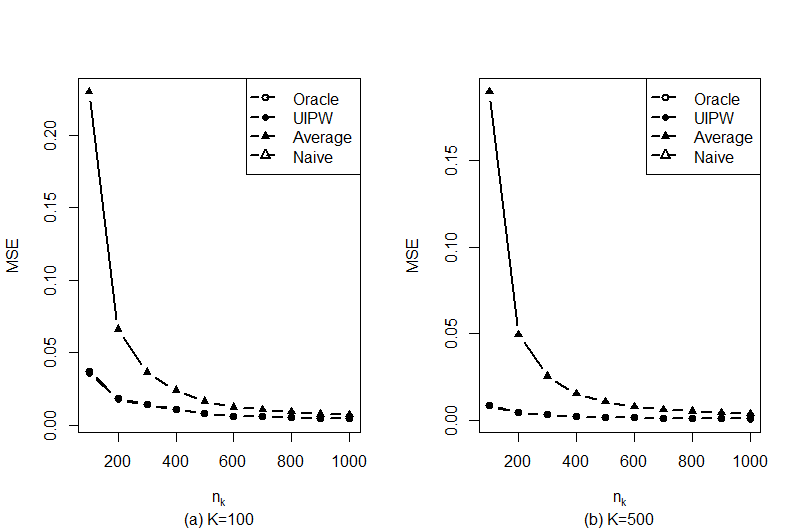}
	}

	\subfigure{
        \centering
		\includegraphics[width=12cm,height=6.5cm]{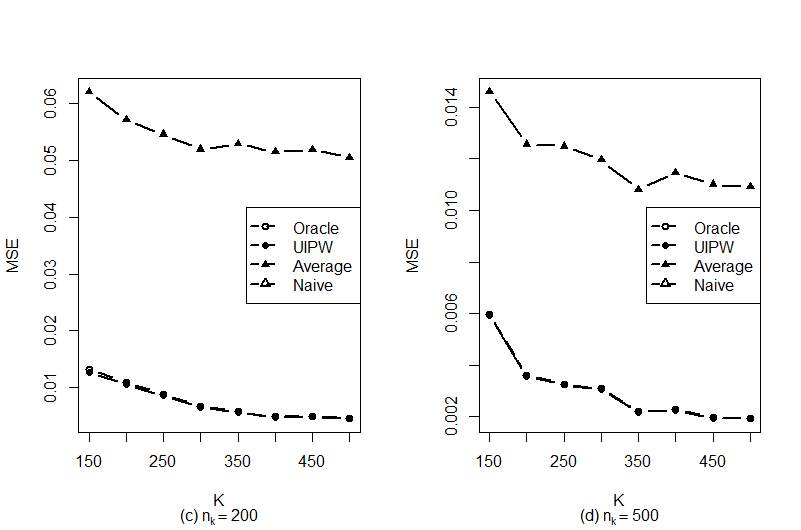}
	}
    \subfigure{
            \centering
			\includegraphics[width=12cm,height=6.5cm]{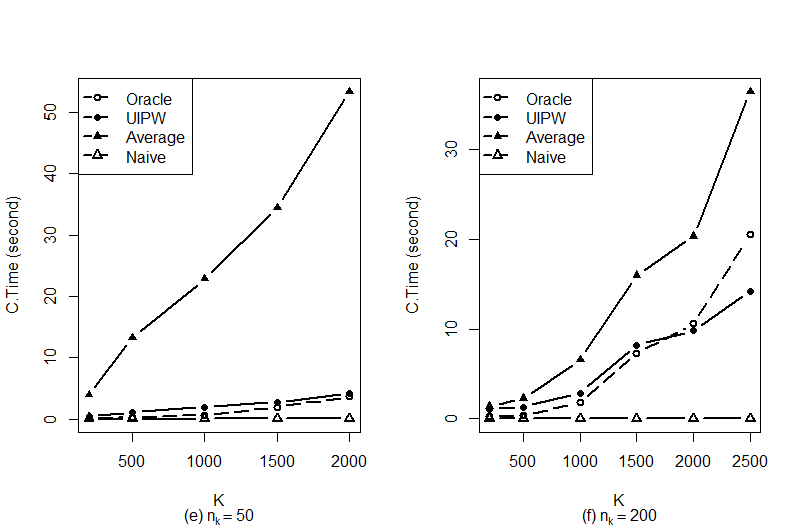}
	}
    \label{logisticMSEp}
\end{figure}

We assume $\bbf d_j=\{\boldsymbol{\delta}_j,\mathbf{Y}_j,\mathbf{X}_j\}$, $j = 1, \ldots, K$ with $\boldsymbol{\delta}_j=(\delta_1^{(j)},\ldots,\delta_{n_j}^{(j)})^{\top}$ is the indicator vector, $\mathbf{Y}_j=(\mathrm{y}_{1}^{(j)},\ldots,\mathrm{y}_{n_j}^{(j)})^{\top}$ is a binary response subject to missing at random and $\mathbf{X}_j=(\mathrm{x}_{1}^{(j)},\ldots,\mathrm{x}_{n_j}^{(j)})^{\top}$ is the covariate vector fully observed, where $\mathrm{y}_{i}^{(j)} \mid \mathrm{x}_{i}^{(j)}$ are independently sampled from a Bernoulli distribution with probability of success $\pi_j(\mathbf{X}_j)=P(\mathbf{Y}_j=\bbf 1\mid \mathbf{X}_j)$.
A logistic regression model takes the form:
\begin{align}
    \pi(\mathbf{X}_j;\boldsymbol{\beta})=\frac{1}{1+e^{-\mathbf{X}_j^{\top}\boldsymbol{\beta}}}. \label{eq:logistic model}
\end{align}
The score function and corresponding negative Hessian matrix for the $j$-th batch of data $\boldsymbol{d}_j$ are respectively written as $\mathbf{S}(\boldsymbol{\beta}|\bbf d_j,\boldsymbol{\alpha})=(1/n_j)\mathbf{X}_j^{\top}\mathbf{D}_j(\boldsymbol{\alpha})(\pi_j(\mathbf{X}_j;\boldsymbol{\beta})-\mathbf{Y}_j)$ and $\mathbf{R}(\boldsymbol{\beta}|\bbf d_j,\boldsymbol{\alpha})=(1/n_j)\mathbf{X}_j^{\top}\mathbf{Q}_j(\boldsymbol{\alpha,\boldsymbol{\beta}})\mathbf{X}_j$, where
\begin{align*}
    \mathbf{D}_j(\boldsymbol{\alpha})&=\mbox{diag}\bigg(\frac{\delta_i}{\pi(\mathrm{x}_{i}^{(j)};\boldsymbol{\alpha})}\bigg),i=1,\ldots,n_j, \\
    \mathbf{Q}_j(\boldsymbol{\alpha},\boldsymbol{\beta})&=\mbox{diag}\bigg(\frac{\delta_i}{\pi(\mathrm{x}_{i}^{(j)};\boldsymbol{\alpha})}\pi_j(\mathbf{X}_j;\boldsymbol{\beta})(1-\pi_j(\mathbf{X}_j;\boldsymbol{\beta}))\bigg),i=1,\ldots,n_j.
\end{align*}
In this example, when the $j$-th batch of data arrives, we consider the following logistic regression model:
\begin{align}
    \pi(\mathbf{X}_j;\boldsymbol{\beta}^0)=\frac{1}{1+exp(-(0.5\mathrm{x}_{i1}^{(j)}+0.5\mathrm{x}_{i2}^{(j)}+\mathrm{x}_{i3}^{(j)}+\mathrm{x}_{i4}^{(j)}))}, \quad i \in \bbf i_j \cap \bbf O,
\end{align}
where $\mathrm{x}_{\cdot t}^{(j)} \sim U(0,1), t=1,2$, $\mathrm{x}_{\cdot t}^{(j)} \sim N(0,1), t=3,4$, and $\boldsymbol{\beta}^0=(0.5,0.5,1,1)$. The propensity function is set to be:
\begin{align*}
    \pi_0(\mathbf{X}_j;\boldsymbol{\alpha}^0)=\frac{1}{1+exp(-(0.5\mathrm{x}_{i1}^{(j)}+\mathrm{x}_{i2}^{(j)}+1.5\mathrm{x}_{i3}^{(j)}+0.5\mathrm{x}_{i4}^{(j)}))},\quad i=1,\ldots,n_j
\end{align*}
The simulation results are reported in Table \ref{logisticMSEt} and Figure \ref{logisticMSEp}.
We have the following findings:
\begin{itemize}
    \item[(1)] \textit{Mean squared error.} Under the criterion of MSE, in Scenario 1, as shown in Table \ref{logisticMSEt}, the oracle estimator and our UIPW estimator are more robust to the varing of $K$ than the simple average IPW estimator.
    In Scenario 2, Figure \ref{logisticMSEp} (a)-(d) implies that the MSE of our UIPW estimator decreases with the increase of $n_k$ or $K$, and moreover, our UIPW estimator is significantly superior to the simple average IPW estimator and has similar behavior to the oracle estimator. 
    \item[(2)] \textit{Computation time.} Under the criterion of C.Time, the naive estimator takes the least C.Time beacuse it only uses the current data.
    In Scenario 1, when $N_k$ is fixed, Table \ref{logisticMSEt} indicates that when the batch size is as small as $n_k=50$, as we said before, the C.Times of the simple average IPW estimator and the naive IPW estimator are larger due to the low convergence rate of the SGD algorithm.
    In Scenario 2, when $n_k$ is small, Figure \ref{logisticMSEp} (e) implies that the simple average estimator takes the most C.Time and our UIPW estimator takes more C.Time than the oracle estimator. 
    However, setting a large $n_k$, as shown in Figure \ref{logisticMSEp} (f), the simple average IPW estimator still takes the most C.Time, but the C.Time of the oracle estimator increases rapidly with the increase of $K$, and gradually exceeds that of our UIPW estimator.
\end{itemize}

\begin{table}[hbp]
    \centering
    \caption{The MSE and computation time are summarized over 200 replication, under the setting of $N_K=100,000$ and $p = 4$ for the logistic regression model (\ref{eq:logistic model}) with batch size $n_k$ from 50 to 2000.}
    \begin{tabular}{|ccccc|}
        \hline
        \multirow{2}*{}         &\multicolumn{4}{c|}{$K=50$, $n_k=2000$}\\
                 &Oracle&UIPW&Average&Naive\\
        MSE      &3.39$\times$$10^{-3}$&3.39$\times$$10^{-3}$&5.53$\times$$10^{-3}$&1.20$\times$$10^{-1}$\\
        C.Time(s)&0.39&20.98&24.12&0.47\\
        \hline
        \multirow{2}*{} &\multicolumn{4}{c|}{$K=100$, $n_k=1000$}\\
        &Oracle&UIPW&Average&Naive\\
        MSE      &3.91$\times$$10^{-3}$&3.89$\times$$10^{-3}$&6.25$\times$$10^{-3}$&4.16$\times$$10^{-1}$\\
        C.Time(s)&0.31&9.43&11.5&0.11\\
        \hline
        \multirow{2}*{}         &\multicolumn{4}{c|}{$K=200$, $n_k=500$}\\
                 &Oracle&UIPW&Average&Naive\\
        MSE      &4.33$\times$$10^{-3}$&4.28$\times$$10^{-3}$&1.31$\times$$10^{-2}$&7.65$\times$$10^{-1}$\\
        C.Time(s)&0.73&5.58&7.98&0.03\\
        \hline
        \multirow{2}*{} &\multicolumn{4}{c|}{$K=500$, $n_k=200$}\\
        &Oracle&UIPW&Average&Naive\\
        MSE      &5.03$\times$$10^{-3}$&4.93$\times$$10^{-3}$&5.14$\times$$10^{-2}$&2.11\\
        C.Time(s)&1.47&3.02&4.88&0.01\\
        \hline
        \multirow{2}*{}         &\multicolumn{4}{c|}{$K=1000$, $n_k=100$}\\
                 &Oracle&UIPW&Average&Naive\\
        MSE      &5.22$\times$$10^{-3}$&5.04$\times$$10^{-3}$&1.9$\times$$10^{-1}$&5.52\\
        C.Time(s)&1.61&2.55&5.52&0.01\\
        \hline
        \multirow{2}*{}&\multicolumn{4}{c|}{$K=2000$, $n_k=50$}\\
        &Oracle&UIPW&Average&Naive\\
        MSE      &1.24$\times$$10^{-2}$&1.27$\times$$10^{-2}$&9.65$\times$$10^{-2}$&3.64\\
        C.Time(s)&3.57&3.83&48.4&0.075\\
        \hline
    \end{tabular}
    \label{logisticMSEt}
\end{table}

\subsubsection{Heteroscedastic logistic regression model}
\begin{figure}[htbp]
	\caption{The MSE in the logistic regression model (\ref{eq:holgistic model}). Note that the MSE of the naive estimator is not depicted in (a) and (b) because it is out of the designed figures.}
	\centering
    \vspace{-0.2cm}
    \subfigbottomskip=-2cm
	\subfigure{
        \centering
        \includegraphics[width=12cm,height=6.5cm]{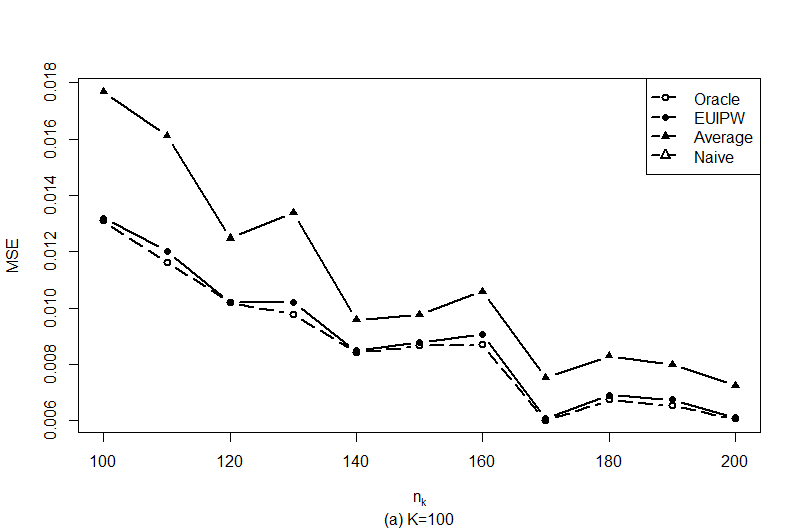}
	}

	\subfigure{
        \centering
		\includegraphics[width=12cm,height=6.5cm]{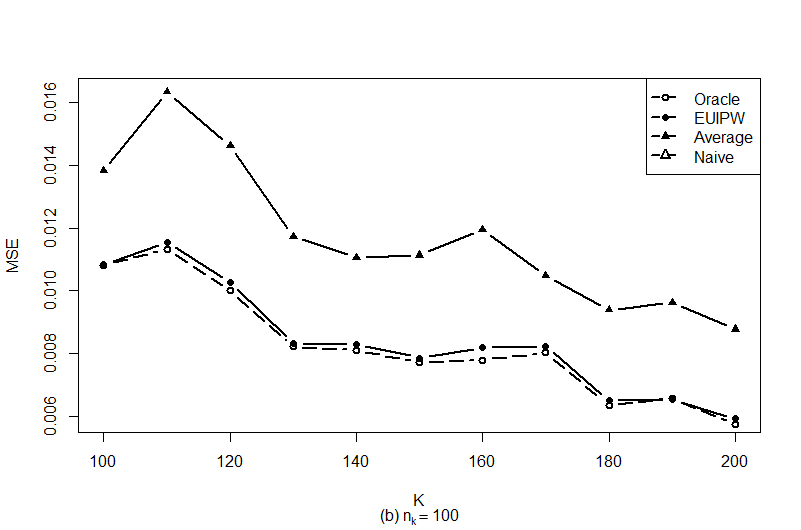}
	}
    \label{hlogisticMSEp}
\end{figure}

We consider a logistic regression between a binary outcome $\mathbf{Y}$ and an exposure $\mathbf{X}$, controlling for a confounding variable $\mathbf{Z}$.
We assume $\bbf d_j=\{\boldsymbol{\delta}_j,\mathbf{Y}_j,\mathbf{X}_j,\mathbf{Z}_j\}$, $j = 1, \ldots, K$ with $\boldsymbol{\delta}_j=(\delta_1^{(j)},\ldots,\delta_{n_j}^{(j)})^{\top}$ is the indicator vector, $\mathbf{Y}_j=(\mathrm{y}_{1}^{(j)},\ldots,\mathrm{y}_{n_j}^{(j)})^{\top}$ is a binary response subject to missing at random, $\mathbf{X}_j=(\mathrm{x}_{1}^{(j)},\ldots,\mathrm{x}_{n_j}^{(j)})^{\top}$ and $\mathbf{Z}_j=(\mathrm{z}_{1}^{(j)},\ldots,\mathrm{z}_{n_j}^{(j)})^{\top}$ are the covariates fully observed,
It is assumed that, for the $j$-th batch of data,
\begin{align}
    \text{logit}\{P(\mathbf{Y}_j=\boldsymbol{1}\mid \mathbf{X}_j, \mathbf{Z}_j)\} = \mathbf{X}_j^{\top}\boldsymbol{\beta} + \mathbf{Z}_j^{\top}\boldsymbol{\gamma}_j, \label{eq:holgistic model}
\end{align}
where $\boldsymbol{\gamma}_j=(\gamma_{1j},\gamma_{2j})^{\top}$ is the nuisance parameter. We set the true value of $\boldsymbol{\beta} = (0.5,0.5)$. The nuisance parameters $\gamma_{1j}$ and $\gamma_{2j}$ are generated
from the uniform distributions $U(-1,1)$ and $U(-2,2)$, respectively.
We generate $\mathbf{X}_j \sim U(0,1)$ and $\mathbf{Z}_j \sim N(\mathbf{X}_j-0.3,1)$. The propensity function is set to be:
\begin{align*}
    \pi_0(\mathbf{X}_j;\boldsymbol{\alpha}^0)=
    \frac{1}{1+\exp(-(0.5\mathrm{x}_{i1}^{(j)}+\mathrm{x}_{i2}^{(j)}))},
\end{align*}
where $\boldsymbol{\alpha}^0=(0.5,1)$.
The simulation results are showed in Table \ref{hlogisticMSEt} and Figure \ref{hlogisticMSEp}. We have the following findings:\
\begin{itemize}
    \item[(1)] Under the criterion of MSE, in Scenario 1, Table \ref{hlogisticMSEt} indicates that the MSE of our EUIPW estimator is close to that of the oracle estimator, while the MSE of the simple average IPW estimator is large. 
    In Scenario 2, when $n_k$ or $K$ increases, Figure \ref{hlogisticMSEp} implies that our EUIPW estimator always exhibits similar performances to the oracle estimator. Both of their MSEs tend to decline, but the decline process is not stable. 
    \item[(2)] Under the criterion of C.Time, as shown in Table \ref{hlogisticMSEt}, the oracle estimator takes the most C.Time because of the computation of the matrix inverse in (\ref{eq:ES}), hence our EUIPW estimator shows clear advantages.
\end{itemize}

\begin{table}[htbp]
    \centering
    \caption{The MSE and computation time are summarized over 200 replication, under the setting of $N_K=20,000$ and $d = 6, p=2, q=2$ for the logistic regression model (\ref{eq:holgistic model}) with batch size $n_k$ from 40 to 100.}
    \begin{tabular}{|ccccc|}
        \hline
        \multirow{2}*{}         &\multicolumn{4}{c|}{$K=40$, $n_k=500$}\\
                 &Oracle&EUIPW&Average&Naive\\
        MSE      &3.98$\times$$10^{-3}$&4.16$\times$$10^{-3}$&4.22$\times$$10^{-3}$&1.08$\times$$10^{-1}$\\
        C.Time(s)&129.99&4.31&2.7&0.09\\
        \hline
        \multirow{2}*{} &\multicolumn{4}{c|}{$K=50$, $n_k=400$}\\
        &Oracle&EUIPW&Average&Naive\\
        MSE      &5.24$\times$$10^{-3}$&5.31$\times$$10^{-3}$&5.72$\times$$10^{-3}$&3.2$\times$$10^{-1}$\\
        C.Time(s)&115.13&3.02&2.92&0.08\\
        \hline
        \multirow{2}*{}         &\multicolumn{4}{c|}{$K=80$, $n_k=250$}\\
                 &Oracle&EUIPW&Average&Naive\\
        MSE      &7.11$\times$$10^{-3}$&7.28$\times$$10^{-3}$&8.33$\times$$10^{-3}$&4.6$\times$$10^{-1}$\\
        C.Time(s)&142.86&1.68&1.51&0.05\\
        \hline
        \multirow{2}*{} &\multicolumn{4}{c|}{$K=100$, $n_k=200$}\\
        &Oracle&EUIPW&Average&Naive\\
        MSE      &8.99$\times$$10^{-3}$&9.08$\times$$10^{-3}$&1.03$\times$$10^{-2}$&1.65\\
        C.Time(s)&117&1.4&1.68&0.08\\
        \hline
    \end{tabular}
    \label{hlogisticMSEt}
\end{table}

\subsection{Real data analysis}
We apply the proposed method UIPW to the National Alzheimer’s Coordinating Center (NACC) Uniform Data Set.
We choose the age (NACCAGE), diabetes (DIABETES, yes/no), depression or dysphoria (DEPD, yes/no) and Mini-Mental State Exam (NACCMMSE) as covariates according to Chen \& Zhou (2011).
As introduced by NACC, we determine whether someone has Alzheimer's disease (AD) based on presumptive etiologic diagnosis of the cognitive disorder (NACCALZD) and cognitive status (NACCUDSD).
That is, someone has AD if NACCALZD = 1 and NACCUDSD = 4, otherwise not.
Let the response $Y = 1$ if someone has AD, otherwise $Y = 0$.
Our goal is to analyze the relationship between selected covariates and AD.

In this example, the streaming data set were formed by monthly visitor data from the period of 7 years over January 2008 to March 2015, with $K=87$ batches of data and a total simple size $N_K=69688$.
However, patients may miss a clinic visit or refuse to undergo a clinical examination during the clinic visit, leading to the incomplete responses.
According to Chen \& Zhou (2011) the MAR mechanism is reasonable.
There are 10359 subjects with missing response accounts for 14.8\%.
Because the response $Y$ is a binary outcome, we utilize (\ref{eq:simulation propensity function}) and (\ref{eq:logistic model}) as the regression models.
We apply our proposed method to construct updatable sequentially parameter estimates.

\begin{figure}[htbp]
	\caption{Trace plots for the coefficient estimates and 95\% pointwise confidence bands of regression coefficients. Note that the confidence bands in (b) is too narrow to show clearly.}
	\centering
    \includegraphics[width=15cm,height=10cm]{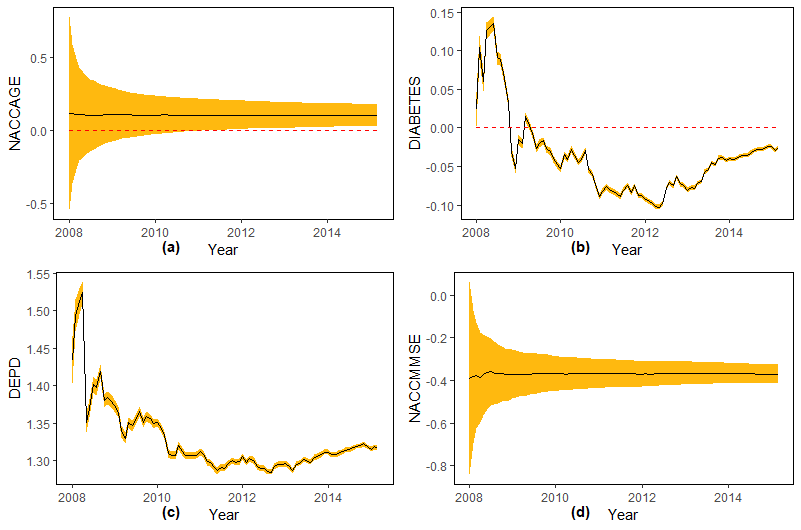}
    \label{confidence_interval}
\end{figure}

Figure \ref{confidence_interval} depicts the 95\% pointwise confidence bands of all regression coefficients.
It is seen that with more batches of data arrived, the confidence bands became narrower.
In Figure \ref{confidence_interval}, (a) and (c) show that AD is positively correlated with NACCAGE and DEPD, meaning that older or more depressed or dysphoric people have a higher prevalence.
It is interesting to find that the trace plot for the DIABETES shows a downward trend from positive to negative as the sample size increases.
This suggests that diabetes now has a positive effect to protect the occurrence of AD. 
Figure \ref{confidence_interval} (d) shows that MMSE has a positive effect to protect the occurrence of AD.

\begin{figure}[htbp]
	\caption{The proportion of correct classifications for the UIPW method and oracle method.}
	\centering
    \includegraphics[width=15cm,height=10cm]{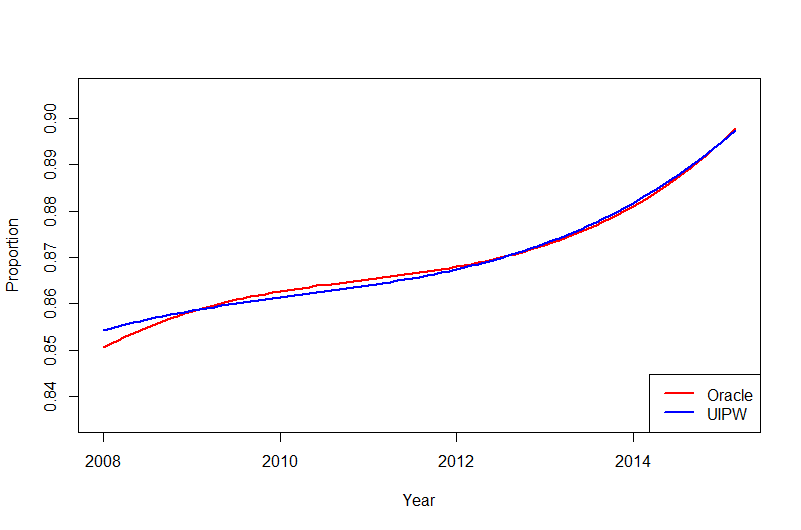}
    \label{Proportion}
\end{figure}

Furthermore, the complete data in the above streaming data set is used to evaluate the UIPW estimator.
Figure \ref{Proportion} shows the proportion of correct classifications for the UIPW method and the oracle method when the classification threshold is 0.5.
It is seen that the proportion of correct classifications for the UIPW method increases when more data are obtained.
To accurately evaluate the performance of the UIPW method and the oracle method as classifiers, we create receiver operating characteristic (ROC) curves using classification probability thresholds between 0 and 1, and then calculate the areas under the ROC curves (AUC).
AUC is one of the most important evaluation metrics for measuring the performance of any classification model.
It is a performance measurement for a classification problem at various thresholds settings.
The larger the value of AUC, the better the effect of classifier.
As shown in Table \ref{AUC}, both the AUCs are between 0.85 and 0.95, meaning that both the UIPW method and the oracle method have good performance.

\begin{table}[htbp]
    \centering
    \caption{The AUCs for the UIPW classifier and the oracle classifier}
        \begin{tabular}{cc}
            \toprule
            Method&AUC\\
            \midrule
            UIPW&0.9167\\
            Oracle&0.9173\\
            \bottomrule
        \end{tabular}
        \label{AUC}
\end{table}

\section{Conclusions and future works}
As shown in Introduction, although a large number of statistical methods and computational recipes have been developed to address the challenges of analyzing the models with streaming data sets or missing data separately, the strategy of online updating IPW estimation with missing response in streaming data sets has not been built in the existing literature.
To address these issues, we propose the updatable inverse probability weighting (UIPW) estimation via a two-step online updating algorithm in the previous sections.
Both the inverse probability weights and the estimator of the parameter of interest can update simultaneously.
This is a unified framework because our two-step online updating algorithm can not only deal with the problem concerned in this paper, but also is useful for all models with two types of parameters, one of which is ragarded as the nuisance parameter. 
Moreover, our UIPW estimator overcomes the memory constraint and is adaptive to the situation where streaming data sets arrive fast and perpetually.
Both the proposed statistical methodology and the computational algorithms have been justified theoretically and examined numerically in the setting of GLMs.
Summary statistics involved in our proposed method have been verified asymptotically equivalent to the classical sufficient statistic.

In order to extend the UIPW estimation under heterogeneous conditions, motivated by theories of efficient score, the efficient score function is constructed and then the efficient updatable inverse probability weighted (EUIPW) estimation is established.
Through the simulation study in Section 4, it is demonstrated that EUIPW estimation is much better than the competitors.

The research can be further improved from the following aspects.
Firstly, this paper mainly focused on the parametric model, for the nonparametric model, we will encounter new challenges such as the curse of dimensionality and the bandwidth selection.
Secondly, it is natural to think of constructing the updatable imputation method.
In the offline framework, most of the imputation methods impute missing response by estimating the conditional expectation of $Y$ given $X$ or conditional density of $Y$ given $X$. 
However, in the online framework, it is difficult to derive the updatable estimator of the conditional expectation of $Y$ given $X$ or conditional density of $Y$ given $X$ when we don't store the historical data. 
Wellenzohn et al. (2017) proposed Top-k Case Matching (TKCM) to impute missing values in streams of time series data.
But TKCM defines for each time series a set of reference time series and exploits similar historical situations in the reference time series for the imputation.
However, we focus on streaming data sets without any reference streaming data sets.
Furthermore, the application of augmented IPW method in streaming data sets can be researched. 
Thirdly, Taylor et al. (2022) researched the transfer learning when the distributions of $k$-th batch of data and $(k+1)$-th batch of data are the same in some aspects of distribution. This is instructive for us to further study UIPW estimation for the heterogeneous scenario.
These are interesting issues and worthy of further investigation in the future.
\section{Acknowledgements}

The research was supported by the National Key R\&D Program of China (grant No. 2018YFA0703900), the National Natural Science Foundation of China (grant No. 11971265) and the National Statistical Science Research Project (grant No. 2022LD03).

\noindent{\bf Appendix A }

\noindent A.1. Proof of Theorem 2

Assume that conditions C1-C7 given in Section 2.2 hold.
When the $k$-th batch of data arrives, the oracle estimator $\hat{\boldsymbol{\beta}}_{N_k}^*$ is the solution of
\begin{align*}
    \sum_{j=1}^k\mathbf{S}(\boldsymbol{\beta}|\bbf d_j,\boldsymbol{\alpha}_{N_k}^*)=\bbf 0,
\end{align*}
where $\boldsymbol{\alpha}_{N_k}^*=\underset{\boldsymbol{\alpha} \in \mathbb{R}^{p}}{\arg \max } \sum_{j=1}^k\sum_{i \in \bbf i_j} \log\bigg\{\pi(X_i;\boldsymbol{\alpha})^{\delta_i}\big[1-\pi(X_i;\boldsymbol{\alpha})\big]^{(1-\delta_i)}\bigg\}$.

Let $\boldsymbol{\beta}^0$ be the true value of the parameter $\boldsymbol{\beta}$ and $\hat{\boldsymbol{\beta}}_k$ be the updatable estimator. For the first batch of data $\bbf d_1$, we have $\hat{\boldsymbol{\beta}}_1=\hat{\boldsymbol{\beta}}_{N_1}^*$.
Next we prove the consistency of $\hat{\boldsymbol{\beta}}_k$ when $k \geq 2$ by the method of induction.

Define a function
\begin{align}\begin{split}
    g_k(\boldsymbol{\beta})=\frac{1}{N_k}\sum_{j=1}^{k-1}\mathbf{L}_{1}(\hat{\boldsymbol{\beta}}_j|\bbf d_j,\hat{\boldsymbol{\alpha}}_j)(\hat{\boldsymbol{\alpha}}_k-\hat{\boldsymbol{\alpha}}_{k-1})+\frac{1}{N_k}\sum_{j=1}^{k-1}\mathbf{L}_{2}(\hat{\boldsymbol{\beta}}_j|\bbf d_j,\hat{\boldsymbol{\alpha}}_j)(\boldsymbol{\beta}-\hat{\boldsymbol{\beta}}_{k-1})\\
    +\frac{1}{N_k}\mathbf{S}(\boldsymbol{\beta}|\bbf d_k,\hat{\boldsymbol{\alpha}}_k).\label{eq:gk}
\end{split}
\end{align}
According to equation (\ref{eq:update beta}), the updatable estimator $\hat{\boldsymbol{\beta}}_k$ satisfies
\begin{align}
    g_k(\hat{\boldsymbol{\beta}}_k)=0. \label{eq:gk=0}
\end{align}

When $\hat{\boldsymbol{\beta}}_{k-1}$ is consistent, we have
\begin{align}\begin{split}
    g_k(\boldsymbol{\beta}^0)=\frac{1}{N_k}\sum_{j=1}^{k-1}\mathbf{L}_{1}(\hat{\boldsymbol{\beta}}_j|\bbf d_j,\hat{\boldsymbol{\alpha}}_j)(\hat{\boldsymbol{\alpha}}_k-\hat{\boldsymbol{\alpha}}_{k-1})+\frac{1}{N_k}\sum_{j=1}^{k-1}\mathbf{L}_{2}(\hat{\boldsymbol{\beta}}_j|\bbf d_j,\hat{\boldsymbol{\alpha}}_j)(\boldsymbol{\beta}^0-\hat{\boldsymbol{\beta}}_{k-1})\\
    +\frac{1}{N_k}\mathbf{S}(\boldsymbol{\beta}^0|\bbf d_k,\hat{\boldsymbol{\alpha}}_k)=o_p(1). \label{eq:gk0}
\end{split}
\end{align}

Taking the difference between equations (\ref{eq:gk0}) and (\ref{eq:gk=0}), we attain
\begin{align}\begin{split}
    g_k(\boldsymbol{\beta}^0)-g_k(\hat{\boldsymbol{\beta}}_k)=\frac{1}{N_k}\sum_{j=1}^{k-1}\mathbf{L}_{2}(\hat{\boldsymbol{\beta}}_j|\bbf d_j,\hat{\boldsymbol{\alpha}}_j)(\boldsymbol{\beta}^0-\hat{\boldsymbol{\beta}}_{k})-\frac{1}{N_k}\mathbf{S}(\hat{\boldsymbol{\beta}}_k|\bbf d_k,\hat{\boldsymbol{\alpha}}_k)\\
    +\frac{1}{N_k}\mathbf{S}(\boldsymbol{\beta}^0|\bbf d_k,\hat{\boldsymbol{\alpha}}_k)=o_p(1). \label{eq:gk0-gk}
\end{split}
\end{align}

Taking the first-order Taylor series expansion of $\mathbf{S}(\hat{\boldsymbol{\beta}}_k|\bbf d_k,\hat{\boldsymbol{\alpha}}_k)$ in equation (\ref{eq:gk0-gk}) around $\boldsymbol{\beta}^0$, we have
\begin{align}\begin{split}
    \mathbf{S}(\hat{\boldsymbol{\beta}}_k|\bbf d_k,\hat{\boldsymbol{\alpha}}_k)=\mathbf{S}(\boldsymbol{\beta}^0|\bbf d_k,\hat{\boldsymbol{\alpha}}_k)+\bigg[\mathbf{R}_{\boldsymbol{\beta}}(\boldsymbol{\beta}^0|\bbf d_k,\hat{\boldsymbol{\alpha}}_k)-\mathbf{R}_{\boldsymbol{\beta}}(\boldsymbol{\beta}^0|\bbf d_k,\hat{\boldsymbol{\alpha}}_k)\\
    +\mathbf{R}_{\boldsymbol{\beta}}(\boldsymbol{\phi}_k|\bbf d_k,\hat{\boldsymbol{\alpha}}_k)\bigg](\hat{\boldsymbol{\beta}}_k-\boldsymbol{\beta}^0),\label{eq:S Taylor expansion}
\end{split}
\end{align}
where  $\boldsymbol{\phi}_k$ lies in between $\hat{\boldsymbol{\beta}}_k$ and $\boldsymbol{\beta}^0$.
By the Lipschitz continuity in condition C7, there exists $\boldsymbol{C}(\boldsymbol{d}_k) \textgreater 0$ such that
\begin{align}
    \parallel\mathbf{R}_{\boldsymbol{\beta}}(\boldsymbol{\phi}_k|\boldsymbol{d}_k,\hat{\boldsymbol{\alpha}}_k,)-\mathbf{R}_{\boldsymbol{\beta}}(\boldsymbol{\beta}^0|\boldsymbol{d}_k,\hat{\boldsymbol{\alpha}}_k,)\parallel \leq \boldsymbol{C}(\boldsymbol{d}_k)\parallel\boldsymbol{\phi}_k-\boldsymbol{\beta}^0\parallel \leq \boldsymbol{C}(\boldsymbol{d}_k)\parallel\hat{\boldsymbol{\beta}}_k-\boldsymbol{\beta}^0\parallel.\label{eq:Lipschitz continuity}
\end{align}
Then we can rewrite equation (\ref{eq:S Taylor expansion}) as
\begin{align}
    \mathbf{S}(\hat{\boldsymbol{\beta}}_k|\bbf d_k,\hat{\boldsymbol{\alpha}}_k)=\mathbf{S}(\boldsymbol{\beta}^0|\bbf d_k,\hat{\boldsymbol{\alpha}}_k)+\mathbf{R}_{\boldsymbol{\beta}}(\boldsymbol{\beta}^0|\bbf d_k,\hat{\boldsymbol{\alpha}}_k)(\hat{\boldsymbol{\beta}}_k-\boldsymbol{\beta}^0)+O_p(n_k\parallel\hat{\boldsymbol{\beta}}_k-\boldsymbol{\beta}^0\parallel^2).\label{eq:simplify Taylor expansion}
\end{align}
Combining equations (\ref{eq:gk0-gk}) and (\ref{eq:simplify Taylor expansion}) yields
\begin{align}\begin{split}
    g_k(\boldsymbol{\beta}^0)-g_k(\hat{\boldsymbol{\beta}}_k)=\frac{1}{N_k}\bigg[\sum_{j=1}^{k-1}\mathbf{L}_{2}(\hat{\boldsymbol{\beta}}_j|\bbf d_j,\hat{\boldsymbol{\alpha}}_j)+\mathbf{R}_{\boldsymbol{\beta}}(\boldsymbol{\beta}^0|\bbf d_k,\hat{\boldsymbol{\alpha}}_k)\bigg](\hat{\boldsymbol{\beta}}_0-\hat{\boldsymbol{\beta}}_k)\\
    +O_p(\frac{n_k}{N_k}\parallel\hat{\boldsymbol{\beta}}_k-\boldsymbol{\beta}^0\parallel^2)=o_p(1).\label{eq:simplify gk0-gk}
\end{split}
\end{align}

We know that $\hat{\boldsymbol{\alpha}}_j$ is consistent for $j=1,\ldots,k$.
Therefore, by conditions C5, $\frac{1}{N_k}[\sum_{j=1}^{k-1}\mathbf{L}_{2}(\hat{\boldsymbol{\beta}}_j|\bbf d_j,\hat{\boldsymbol{\alpha}}_j)+\mathbf{R}_{\boldsymbol{\beta}}(\boldsymbol{\beta}^0|\bbf d_k,\hat{\boldsymbol{\alpha}}_k)]$ is positive definite. Under the assumption  $\hat{\boldsymbol{\beta}}_j$ is consistent for $j=1,\ldots,k-1$.
It follows that $\hat{\boldsymbol{\beta}}_k-\boldsymbol{\beta}^0\overset{P}{\rightarrow}0$, as $N_k\rightarrow\infty$.

\noindent A.2. Proof of Theorem 3

(a) When $k=1$, that is, for the first batch of data $\bbf d_1$, $n_1=N_1$, the oracle estimator $\hat{\boldsymbol{\beta}}_{N_1}^{*}=\hat{\boldsymbol{\beta}}_1$ satisfies $(1/N_1)\mathbf{S}(\hat{\boldsymbol{\beta}}_1|\boldsymbol{d}_1,\hat{\boldsymbol{\alpha}}_1)=0$ and $\sqrt{N_1}(\hat{\boldsymbol{\beta}}_1-\boldsymbol{\beta}^0)\overset{d}{\rightarrow}N(0,\boldsymbol{\Sigma}^0)$, as $N_1=n_1 \rightarrow \infty$.
Taking the first-order Taylor series expansion of $\mathbf{S}(\boldsymbol{\beta}^0|\bbf d_1,\hat{\boldsymbol{\alpha}}_1)$ around $\hat{\boldsymbol{\beta}}_1$, we have
\begin{align}
    \frac{1}{N_1}\mathbf{S}(\boldsymbol{\beta}^0|\boldsymbol{d}_1,\hat{\boldsymbol{\alpha}}_1)=\frac{1}{N_1}\mathbf{R}_{\boldsymbol{\beta}}(\hat{\boldsymbol{\beta}}_1|\boldsymbol{d}_1,\hat{\boldsymbol{\alpha}}_1)(\boldsymbol{\beta}^0-\hat{\boldsymbol{\beta}}_1)+O_p(\frac{n_1}{N_1}\parallel\hat{\boldsymbol{\beta}}_1-\boldsymbol{\beta}^0\parallel^2). \label{eq:S1 Taylor expansion}
\end{align}
(b)When $k \geq 2$, considering updating $\hat{\boldsymbol{\beta}}_{k-1}$ to $\hat{\boldsymbol{\beta}}_{k}$. The oracle estimator $\hat{\boldsymbol{\beta}}_{N_k}^*$ satisfies $(1/N_k)\sum_{j=1}^k\mathbf{S}(\hat{\boldsymbol{\beta}}_{N_k}^*|\boldsymbol{d}_j,\hat{\boldsymbol{\alpha}}_{N_k}^*)=0$.
According to (\ref{eq:simplify Taylor expansion})
\begin{align}\begin{split}
    \frac{1}{N_k}\sum_{j=1}^k\mathbf{S}(\boldsymbol{\beta}^0|\boldsymbol{d}_j,\hat{\boldsymbol{\alpha}}_{N_k}^*)+\frac{1}{N_k}\sum_{j=1}^k\mathbf{R}_{\boldsymbol{\beta}}(\boldsymbol{\beta}^0|\boldsymbol{d}_j,\hat{\boldsymbol{\alpha}}_{N_k}^*)(\hat{\boldsymbol{\beta}}_{N_k}^*-\boldsymbol{\beta}^0)\\
    +O_p(\parallel\hat{\boldsymbol{\beta}}_{N_k}^*-\boldsymbol{\beta}^0\parallel^2)=\boldsymbol{0}. \label{eq:Sk Taylor expansion}
\end{split}
\end{align}

According to equations (\ref{eq:gk}), (\ref{eq:gk=0}) and (\ref{eq:simplify gk0-gk}), we know that
\begin{align*}
    g_k(\boldsymbol{\beta}^0) =& \frac{1}{N_k}\sum_{j=1}^{k-1}\mathbf{L}_{1}(\hat{\boldsymbol{\beta}}_j|\bbf d_j,\hat{\boldsymbol{\alpha}}_j)(\hat{\boldsymbol{\alpha}}_k-\hat{\boldsymbol{\alpha}}_{k-1})+\frac{1}{N_k}\sum_{j=1}^{k-1}\mathbf{L}_{2}(\hat{\boldsymbol{\beta}}_j|\bbf d_j,\hat{\boldsymbol{\alpha}}_j)(\boldsymbol{\beta}^0-\hat{\boldsymbol{\beta}}_{k-1})\\
    &+\frac{1}{N_k}\mathbf{S}(\boldsymbol{\beta}^0|\bbf d_k,\hat{\boldsymbol{\alpha}}_k) \\
    =& \frac{1}{N_k}\bigg[\sum_{j=1}^{k-1}\mathbf{L}_{2}(\hat{\boldsymbol{\beta}}_j|\bbf d_j,\hat{\boldsymbol{\alpha}}_j)+\mathbf{R}_{\boldsymbol{\beta}}(\boldsymbol{\beta}^0|\bbf d_k,\hat{\boldsymbol{\alpha}}_k)\bigg](\hat{\boldsymbol{\beta}}_0-\hat{\boldsymbol{\beta}}_k)\\
    &+O_p(\frac{n_k}{N_k}\parallel\hat{\boldsymbol{\beta}}_k-\boldsymbol{\beta}^0\parallel^2)=o_p(1).\\
\end{align*}
It follows that
\begin{align}\begin{split}&
    \frac{1}{N_k}\bigg[\sum_{j=1}^{k-1}\mathbf{L}_{2}(\hat{\boldsymbol{\beta}}_j|\bbf d_j,\hat{\boldsymbol{\alpha}}_j)+\mathbf{R}_{\boldsymbol{\beta}}(\boldsymbol{\beta}^0|\bbf d_k,\hat{\boldsymbol{\alpha}}_k)\bigg](\hat{\boldsymbol{\beta}}_k-\boldsymbol{\beta}^0)+\frac{1}{N_k}\sum_{j=1}^{k-1}\mathbf{L}_{2}(\hat{\boldsymbol{\beta}}_j|\bbf d_j,\hat{\boldsymbol{\alpha}}_j)\\
    &\times(\boldsymbol{\beta}^0-\hat{\boldsymbol{\beta}}_{k-1})+\frac{1}{N_k}\mathbf{S}(\boldsymbol{\beta}^0|\bbf d_k,\hat{\boldsymbol{\alpha}}_k)+O_p(\frac{n_k}{N_k}\parallel\hat{\boldsymbol{\beta}}_k-\boldsymbol{\beta}^0\parallel^2)+o_p(1)=\bbf 0. \label{eq:gk0=gk0}
\end{split}
\end{align}
Similarly to equation (\ref{eq:S1 Taylor expansion}), at the $(k-1)$-th batch of data, it is easy to show that
\begin{align}\begin{split}
    \frac{1}{N_{k-1}}\sum_{j=1}^{k-1}\mathbf{S}(\boldsymbol{\beta}^0|\boldsymbol{d}_j,\hat{\boldsymbol{\alpha}}_{k-1})&=\frac{1}{N_{k-1}}\sum_{j=1}^{k-1}\mathbf{S}(\boldsymbol{\beta}^0|\boldsymbol{d}_j,\hat{\boldsymbol{\alpha}}_{k})+o_p(1)
\\&
    =\frac{1}{N_{k-1}}\sum_{j=1}^{k-1}\mathbf{L}_2(\hat{\boldsymbol{\beta}}_j|\boldsymbol{d}_j,\hat{\boldsymbol{\alpha}}_j)(\boldsymbol{\beta}^0-\hat{\boldsymbol{\beta}}_{k-1})
\\&
    +O_p(\sum_{j=1}^{k-1}\frac{n_j}{N_{k-1}}\parallel\hat{\boldsymbol{\beta}}_j-\boldsymbol{\beta}^0\parallel^2)+o_p(1). \label{eq:simplify Sk Taylor expansion}
\end{split}
\end{align}
Plugging equation (\ref{eq:simplify Sk Taylor expansion}) into equation (\ref{eq:gk0=gk0}), we obtain
\begin{align}\begin{split}&
    \frac{1}{N_{k}}\sum_{j=1}^{k}\mathbf{S}(\boldsymbol{\beta}^0|\boldsymbol{d}_j,\hat{\boldsymbol{\alpha}}_k)+\frac{1}{N_k}\bigg\{\sum_{j=1}^{k-1}\mathbf{L}_2(\hat{\boldsymbol{\beta}}_j|\boldsymbol{d}_j,\hat{\boldsymbol{\alpha}}_j)+\mathbf{R}_{\boldsymbol{\beta}}(\boldsymbol{\beta}^0|\boldsymbol{d}_k,\hat{\boldsymbol{\alpha}}_k)\bigg\}(\hat{\boldsymbol{\beta}}_{k}-\boldsymbol{\beta}^0)
    \\&
    +\frac{1}{N_k}\sum_{j=1}^{k-1}\mathbf{L}_{1}(\hat{\boldsymbol{\beta}}_j|\bbf d_j,\hat{\boldsymbol{\alpha}}_j)(\hat{\boldsymbol{\alpha}}_k-\hat{\boldsymbol{\alpha}}_{k-1})+O_p(\sum_{j=1}^{k}\frac{n_j}{N_{k}}\parallel\hat{\boldsymbol{\beta}}_j-\boldsymbol{\beta}^0\parallel^2)=\boldsymbol{0}.
\end{split}
\end{align}
According to Lemma 1 and Theorem 2, $\hat{\boldsymbol{\alpha}}_j$ are consistent for $j=1,\ldots,k$ and $\hat{\boldsymbol{\beta}}_j$ are consistent for $j=1,\ldots,k-1$.
Then, by condition C7, the continuous mapping theorem implies that
\begin{align}\begin{split}
    \frac{1}{N_{k}}\sum_{j=1}^{k}\mathbf{S}(\boldsymbol{\beta}^0|\boldsymbol{d}_j,\hat{\boldsymbol{\alpha}}_k)+\frac{1}{N_k}\sum_{j=1}^{k}\mathbf{R}_{\boldsymbol{\beta}}(\boldsymbol{\beta}^0|\boldsymbol{d}_j,\hat{\boldsymbol{\alpha}}_j)(\hat{\boldsymbol{\beta}}_{k}-\boldsymbol{\beta}^0)\\
    +O_p(\sum_{j=1}^{k}\frac{n_j}{N_{k}}\parallel\hat{\boldsymbol{\beta}}_j-\boldsymbol{\beta}^0\parallel^2) = \boldsymbol{0}. \label{eq:CM SK}
\end{split}
\end{align}
By condition C5, $\mathcal{I}^{-1}_{N_{k}}(\boldsymbol{\beta}^0)$ exists, and thus the central limit theorem implies that
\begin{align}
    \sqrt{N_k}(\hat{\boldsymbol{\beta}}_k-\boldsymbol{\beta}^0)=\{-\frac{1}{N_k}\sum_{j=1}^k\mathbf{R}_{\boldsymbol{\beta}}(\boldsymbol{\beta}^0|\boldsymbol{d}_j,\hat{\boldsymbol{\alpha}}_j)\}^{-1}\frac{1}{\sqrt{N_k}}\sum_{j=1}^k\mathbf{S}(\boldsymbol{\beta}^0|\boldsymbol{d}_j,\hat{\boldsymbol{\alpha}}_k)+o_p(1).\label{eq:beta asymptotic expression 1}
\end{align}
Taking the first-order Taylor series expansion of $\mathbf{S}(\boldsymbol{\beta}^0|\boldsymbol{d}_j,\hat{\boldsymbol{\alpha}}_k)$ around $\boldsymbol{\alpha}^0$ in (\ref{eq:beta asymptotic expression 1}), we obtain
\begin{align}\begin{split}
    \sqrt{N_k}(\hat{\boldsymbol{\beta}}_k-\boldsymbol{\beta}^0)=\bigg\{-\frac{1}{N_k}\sum_{j=1}^k\mathbf{R}_{\boldsymbol{\beta}}(\boldsymbol{\beta}^0|\boldsymbol{d}_j,\hat{\boldsymbol{\alpha}}_j)\bigg\}^{-1}\bigg\{\frac{1}{\sqrt{N_k}}\sum_{j=1}^k\mathbf{S}(\boldsymbol{\beta}^0|\boldsymbol{d}_j,\boldsymbol{\alpha}^0)\\
    +\big[\frac{1}{N_k}\sum_{j=1}^{k}\mathbf{R}_{\boldsymbol{\alpha}}(\boldsymbol{\beta}^0|\boldsymbol{d}_j,\boldsymbol{\psi}_k)\sqrt{N_k}(\hat{\boldsymbol{\alpha}}_k-\boldsymbol{\alpha}^0)\big]\bigg\}+o_p(1). \label{eq:beta asymptotic expression 2}
\end{split}
\end{align}
where  $\boldsymbol{\psi}_k$ lies in between $\hat{\boldsymbol{\alpha}}_k$ and $\boldsymbol{\alpha}^0$. According to the asymptotic normality of $\hat{\boldsymbol{\alpha}}_k$, we can get
\begin{align}
    \sqrt{N_k}(\hat{\boldsymbol{\alpha}}_k-\boldsymbol{\alpha}^0) = \frac{1}{\sqrt{N_k}}\sum_{j = 1}^{k}\bigg\{E\big[V(\delta;X,\boldsymbol{\alpha}^0)V(\delta;X,\boldsymbol{\alpha}^0)^{\top}\big]\bigg\}^{-1}\mathbf{V}(\boldsymbol{\delta}_j;\mathbf{X}_j,\boldsymbol{\alpha}^0)+o_p(1) \label{eq:alpha asymptotic expression 1}
\end{align}
Substituting (\ref{eq:alpha asymptotic expression 1}) into (\ref{eq:beta asymptotic expression 2}) gives
\begin{align}\begin{split}&
    \sqrt{N_k}(\hat{\boldsymbol{\beta}}_k-\boldsymbol{\beta}_0)=\bigg\{-\frac{1}{N_k}\sum_{j=1}^k\mathbf{R}_{\boldsymbol{\beta}}(\boldsymbol{\beta}^0|\boldsymbol{d}_j,\hat{\boldsymbol{\alpha}}_j)\bigg\}^{-1}\frac{1}{\sqrt{N_k}}\sum_{j=1}^k\bigg\{\mathbf{S}(\boldsymbol{\beta}^0|\boldsymbol{d}_j,\boldsymbol{\alpha}^0)\\
    &+E\big[R_{\boldsymbol{\alpha}}(\boldsymbol{\beta}^0|W,\boldsymbol{\alpha}^0)\big] \bigg\{E\big[V(\delta;X,\boldsymbol{\alpha}^0)V(\delta;X,\boldsymbol{\alpha}^0)^{\top}\big]\bigg\}^{-1}\mathbf{V}(\boldsymbol{\delta}_j;\mathbf{X}_j,\boldsymbol{\alpha}^0)\bigg\}+o_p(1),\label{eq:asymptotic distribution}
\end{split}
\end{align}
where $R_{\boldsymbol{\alpha}}(\boldsymbol{\beta}^0|W,\boldsymbol{\alpha}^0)$ is the derivative of $S(\boldsymbol{\beta}^0|W,\boldsymbol{\alpha}^0)$ with respect to $\boldsymbol{\alpha}$.
In addition, for the observed data $W$, we have
\begin{align}
    E\big[S(\boldsymbol{\beta}|W,\boldsymbol{\alpha})\big]=\int S(\boldsymbol{\beta}|W,\boldsymbol{\alpha})p(W;\boldsymbol{\beta},\boldsymbol{\alpha})dW=0, \label{eq:ES=0}
\end{align}
where $p(W;\boldsymbol{\beta},\boldsymbol{\alpha})$ is the joint distribution function.
Taking the derivative of equation (\ref{eq:ES=0})
\begin{align}\begin{split}
    \frac{\partial}{\partial \boldsymbol{\alpha}}E[S(\boldsymbol{\beta}|W,\boldsymbol{\alpha})]&=\int \frac{\partial}{\partial \boldsymbol{\alpha}}S(\boldsymbol{\beta}|W,\boldsymbol{\alpha})p(W;\boldsymbol{\beta},\boldsymbol{\alpha})dW+\int S(\boldsymbol{\beta}|W,\boldsymbol{\alpha})\frac{\partial}{\partial \boldsymbol{\alpha}}p(W;\boldsymbol{\beta},\boldsymbol{\alpha})dW\\
    &=E[R_{\boldsymbol{\alpha}}(\boldsymbol{\beta}|W,\boldsymbol{\alpha})]+\int S(\boldsymbol{\beta}|W,\boldsymbol{\alpha})\frac{\frac{\partial}{\partial \boldsymbol{\alpha}}p(W;\boldsymbol{\beta},\boldsymbol{\alpha})}{p(W;\boldsymbol{\beta},\boldsymbol{\alpha})}p(W;\boldsymbol{\beta},\boldsymbol{\alpha})dW\\
    &=E[R_{\boldsymbol{\alpha}}(\boldsymbol{\beta}|W,\boldsymbol{\alpha})]+E[S(\boldsymbol{\beta}|W,\boldsymbol{\alpha})V(\delta;X,\boldsymbol{\alpha})]=0. \label{eq:partial ES}
\end{split}
\end{align}
By the weak law of large numbers and the consistency of $\hat{\boldsymbol{\alpha}}_j$ for $j=1,\ldots,k$, we have $(1/N_k)\sum_{j=1}^K\mathbf{R}_{\boldsymbol{\beta}}(\boldsymbol{\beta}^0|\bbf d_j,\hat{\boldsymbol{\alpha}}_j)\overset{P}{\rightarrow}E\left\{R_{\boldsymbol{\beta}}(\boldsymbol{\beta}|W,\boldsymbol{\alpha})\right\}$, where $R_{\boldsymbol{\beta}}(\boldsymbol{\beta}^0|W,\boldsymbol{\alpha}^0)=\nabla_{\boldsymbol{\beta}}S(\boldsymbol{\beta}^0|W,\boldsymbol{\alpha}^0)$ is the derivative of $S(\boldsymbol{\beta}^0|W,\boldsymbol{\alpha}^0)$ with respect to $\boldsymbol{\beta}$.
Combining equations (\ref{eq:asymptotic distribution}) and (\ref{eq:partial ES}), the central limit theorem implies that
\begin{align*}\begin{split}
        \sqrt{N_k}(\hat{\boldsymbol{\beta}}_k-\boldsymbol{\beta}^0)\overset{d}{\rightarrow} N\bigg(0,\bigg\{E\big[R_{\boldsymbol{\beta}}(\boldsymbol{\beta}^0|W,\boldsymbol{\alpha}^0)\big]\bigg\}^{-1}Var\bigg\{S(\boldsymbol{\beta}^0|W,\boldsymbol{\alpha}^0)-J(\boldsymbol{\beta}^0,\boldsymbol{\alpha}^0)\bigg\}\\
        \times\bigg\{E\big[R_{\boldsymbol{\beta}}(\boldsymbol{\beta}^0|W,\boldsymbol{\alpha}^0)\big]\bigg\}^{-1^{\top}}\bigg),\text{ as } N_k \rightarrow \infty.
\end{split}
\end{align*}
It is easy to get that
$$
E\left[R_{\boldsymbol{\beta}}(\boldsymbol{\beta}^0|W,\boldsymbol{\alpha}^0)\right]=E\left[X\left(v(\mu)g_{\mu}^2(\mu)\right)^{-1}X^{\top}\right].
$$

\end{document}